\documentclass[12pt, a4paper]{amsart}

\usepackage{amsmath, amsthm, amssymb}
\usepackage{graphicx}
\usepackage{enumerate}
\usepackage[top=2.5cm, bottom=2.5cm, left=2.5cm, right=2.5cm]{geometry}

\newcommand{\M}{\mathbb M}

\newcommand{\x}{\mathbf x}
\newcommand{\y}{\mathbf y}

\newtheorem{definition}{Definition}

\newtheorem{proposition}{\bf Proposition}
\newtheorem{theorem}{\bf Theorem}
\newtheorem{corollary}{\bf Corollary}

\renewenvironment{proof}{\noindent {\bf Proof: }}{\rm\\}
\theoremstyle{definition}
\newtheorem{remark}{Remark}{\rm}
\newtheorem{example}{Example}{\rm}

\usepackage{color}
\usepackage{stmaryrd}

\usepackage[T1]{fontenc}
\usepackage[export]{adjustbox}

\usepackage{bbm}

\usepackage{mathrsfs}
\usepackage{mathtools}

\newcommand{\vertiii}[1]{{\left\vert\kern-0.25ex\left\vert\kern-0.25ex\left\vert #1 
    \right\vert\kern-0.25ex\right\vert\kern-0.25ex\right\vert}}

    \usepackage{etoolbox}
    
    \usepackage{arydshln}

\DeclarePairedDelimiterX{\normi}[1]
  {|\!|\!|}
  {|\!|\!|}
  {\ifblank{#1}{\:\cdot\:}{#1}}

\usepackage{algorithm,algpseudocode}
\usepackage{caption}
\DeclareCaptionLabelFormat{algnonumber}{Algorithm}
\begin{document}

\title[Mixed peak-gain/$H_\infty$ synthesis]{Mixed $L_1/H_\infty$-synthesis for $L_\infty$-stability$^*$}
\author{Pierre APKARIAN$^1$}
\author{Dominikus NOLL$^{2}$}
\thanks{$^1$ONERA, The French Aerospace Lab, Toulouse, France}
\thanks{$^2$Institut de Math\'ematiques, Universit\'e de Toulouse, France\\\hspace*{.4cm}${}^*$International Journal of Robust and Nonlinear Control}

\begin{abstract}
We consider stabilization and performance optimization 
of non-linear controlled systems, where the non-linearity satisfies a sector constraint
asymptotically.
This leads to optimization of the closed loop peak-to-peak 
system norm subject to $H_\infty$-performance constraints. Non-linear controlled systems tuned successfully
by this novel approach are locally exponentially stable and globally BIBO-stable.
\\

\noindent
{\sc Key Words.} BIBO-stability $\cdot$ peak-to-peak norm $\cdot$ asymptotic  constraint $\cdot$ 
boundary feedback control $\cdot$ wave equation $\cdot$ sector non-linearity $\cdot$ IQC
\end{abstract}

\maketitle

\section{Introduction}
The peak-gain, or peak-to-peak norm of a BIBO-stable linear time-invariant system $G$, is the time-domain
$L_\infty$ operator norm,
\begin{equation}
    \label{peak_to_peak}
\|G\|_{\rm pk\_gn} = \sup\{ |G\ast w|_\infty: |w|_\infty \leq 1\},
\end{equation}
where the signal norm
on $L_\infty([0,\infty),\mathbb R^n)$ is $|x|_\infty=\sup_{t\geq 0} \max_{i=1,\dots,n} |x_i(t)|$.
In the SISO case it is also known as the system $L_1$-norm.
As opposed to the more standard $H_2$- or $H_\infty$-norms,
computation or optimization of $\|G\|_{\rm pk\_gn}$ has found only mild attention in the control literature,
even though its importance e.g. for the rejection of persistent
perturbations was recognized \cite{dahleh,dahleh_pears,dahleh_voulg,diaz,linnemann,rutland,sznaier}.
One of the reasons of this disesteem  is probably the link of $\|\cdot\|_{\rm pk\_gn}$
with the $H_\infty$-norm $\|\cdot\|_\infty$, where in the chain
\begin{equation}
\label{chain}
m^{-1/2}\|G\|_\infty \leq  \,\|G\|_{\rm pk\_gn} \leq (2n+1){p}^{1/2} \|G\|_\infty
\end{equation}
the right-hand estimate holds
for real-rational systems $G$ with $n$ poles and $p$ inputs, while the left-hand estimate
is valid even for infinite dimensional well-posed  BIBO-stable systems with $m$ outputs.
This may have been interpreted in the sense
that optimizing $\|G\|_{\rm pk\_gn}$
offers nothing substantial over optimizing $\|G\|_\infty$. In the present work we
show that optimizing $\|G\|_{\rm pk\_gn}$ has genuine scope.

For the purpose of motivation,  we consider a possibly infinite-dimensional Lur'e system, 
where a tunable LTI (Linear Time-Invariant) block is in loop with a sector non-linearity.  By the Small Gain theorem
closed-loop $L_2$-stability  is assured if one succeeds in tuning the LTI-block to satisfy a suitable $H_\infty$-norm or frequency shape constraint.
However, this sufficient conditions may be difficult, or even impossible, to
achieve if the sector is too large. Here our new approach applies and replaces the large
sector by a smaller one, which the non-linearity satisfies only  asymptotically.  Application of a small gain argument now requires working with  the 
time-domain $L_\infty$-norm
instead of the $L_2$-norm. In consequence, the LTI-block is now tuned to
satisfy a constraint in the peak-to-peak norm (\ref{peak_to_peak}). If successful, the non-linear closed loop is 
BIBO stable. Due to the smaller primal sector, this is often easier to achieve
than the original $H_\infty$-constraint, and it is one of the few remaining options for non-linear systems with different attraction regimes.

This approach via {\it asymptotic sectors} may be combined with $H_\infty$-methods to guarantee
local exponential stability along with global BIBO-stability. This leads to a novel type of
mixed peak-gain$/H_\infty$-optimization program.

In order to demonstrate the potential of our method, we discuss feedback control of a wave equation
with a non-linear anti-damping boundary causing instability. This model has been used to control
slipstick vibrations in drilling systems \cite{besselink,challamel:00,roman,saldivar:13,saldivar:16}. Our method allows to prove local exponential stability
in tandem with global BIBO-stability for scenarios, where this was previously impossible, the challenge being to achieve this with finite-dimensional
controllers of simple implementable structure. The second part of the paper extends the concept of
asymptotic constraints to MIMO non-linearities,  highlighting that applications are not limited to the SISO case.

The organization %
is as follows. In Section \ref{sect_lure} we discuss
the case of a sector non-linearity.  An algorithm based on mixed $H_\infty/H_\infty$- and 
peak-gain/$H_\infty$-programs is presented in Section \ref{algorithm}. 
In
Section \ref{best} we show how the aperture of the asymptotic sector may be optimized, a feature
which is not possible with standard sectors.
Section \ref{sect_slipstick} 
discusses the application to the control of slipstick vibrations.
Section \ref{mimo}  resumes theory and extends the asymptotic concept to MIMO non-linearity along with illustrations and applications.
Properties of the peak-to-peak norm and implementation of the mixed programs are discussed in Section \ref{peak}.

\section{Mixed program for a  Lur'e systems}
\label{sect_lure}
For the purpose of motivation
we consider a controlled Lur'e system
with state $x$, control input $u$, measured output $y$, disturbances $w$, and regulated outputs $z$:
\begin{align}
\label{NL1}
G_{nl}:\qquad
\begin{split}
\dot{x}& = Ax + B_pp \,  + B_ww\, + B_uu\\
q&= C_qx \qquad\qquad\qquad \; +D_qu\\
p&= \phi(t,q)\\
z&= C_zx  \quad\; \;\;\;\;\, + D_{zw}w + D_{zu}u \\
y &= C_yx
\end{split}
\end{align}
where $p(t) = \phi(t,q(t))$ is a non-linearity
satisfying $\phi(t,0)=0$, $\frac{\partial \phi}{\partial q}(t,0)=0$ and a
sector constraint
\[
(\phi(t,q)-aq)\cdot(\phi(t,q)-bq)\leq 0
\]
for all $t\geq 0$ and all $q$,
abbreviated $\phi\in {\bf sect}(a,b)$. Since
$a \leq \frac{\partial \phi}{\partial q}(t,0)=0 \leq b$, 
the linearized system is 
\begin{align}
\label{L1}
G: \qquad
\begin{split}
\dot{x}& = Ax \,\, + B_ww\,\, + B_uu\\
z&= C_zx + D_{zw}w+D_{zu}u \\
y &= C_yx
\end{split}
\end{align}
In nominal $H_\infty$-synthesis, we might  interpret the non-linearity as a mere disturbance
and optimize a suitable closed-loop performance channel
$\|T_{zw}(G,K)\|_\infty$ over a class $K\in \mathscr K$ of structured controllers \cite{AN06a},  
with optimal $H_\infty$-controller $K^*\in \mathscr K$ and gain
$\gamma_\infty = \|T_{zw}(G,K^*)\|_\infty$. 

\begin{figure}[ht!]
\includegraphics[scale=0.9]{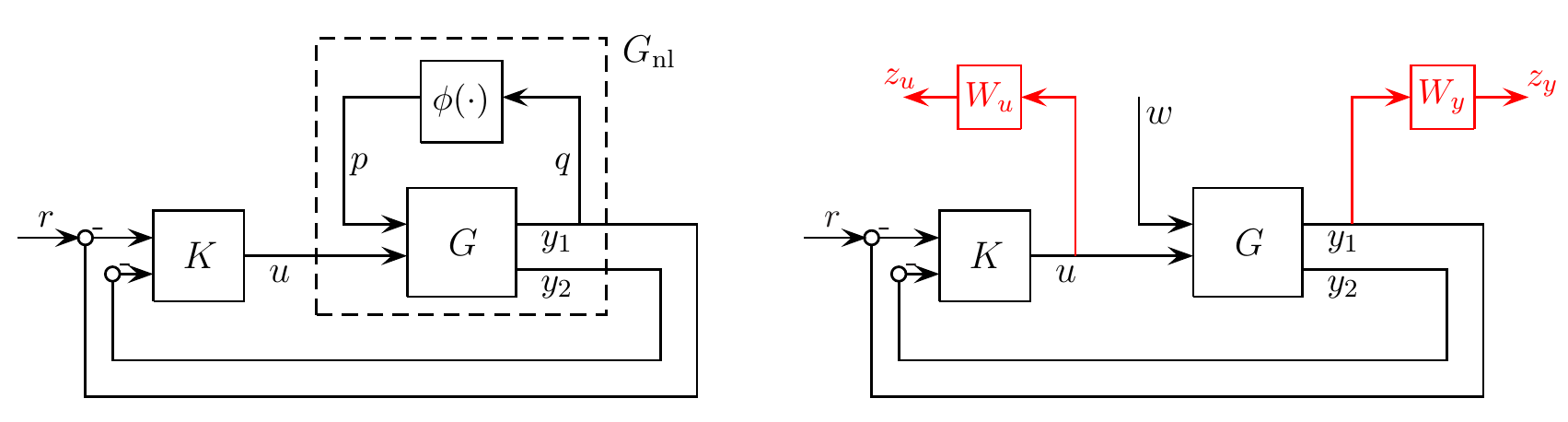}
\caption{Lur'e system $G_{\rm nl}$ as loop between non-linearity $\phi$ and LTI system $G$ (left). Nominal $H_\infty$-synthesis
with $G$ and non-linearity interpreted as disturbance $w$ (right). \label{optimistic}}
\end{figure}

Suppose this optimistic approach of representing the non-linearity by a
disturbance $w$  (as in Fig. \ref{optimistic} right) is too unspecific and $K^*$ is not entirely satisfactory. 
Then we have to target the sector non-linearity explicitly (as in Fig. \ref{optimistic} left). 
Putting $c=(b+a)/2$, $r=(b-a)/2$, and $\psi(t,q) = \phi(t,q)-c q$, we have
$\psi\in {\bf sect}(-r,r)$. The non-linear system (\ref{NL1}) is  now equivalently written as
\begin{align}
\label{NL2}
G_{nl}: \qquad
\begin{split}
\dot{x}& = Ax + cB_p C_q x + B_p p  + B_uu\\
q&= C_qx \qquad \qquad \qquad\;\;\;\, + D_qu \\
p&=\psi(t,q) \\
y &= C_yx
\end{split}
\end{align}
where the performance channel $w\to z$ is temporarily ignored for notational convenience. 
We introduce $A_\psi = A+cB_p C_q$ and %
\begin{align}
\label{L2}
G_{\psi}: \qquad
\begin{split}
\dot{x}& = A_\psi x + B_p p +  B_uu\\
q &= C_q x \qquad\quad  +D_q u\\
y &= C_yx
\end{split}
\end{align}
then (\ref{NL1}) is equivalent to putting
$p\to q$ of $G_{\psi}$ in loop with the centered non-linearity $\psi(\cdot)=\phi(\cdot)-cI$.

Closing the controller loop
$u=Ky$ in $G_{\psi}$ leads to the channel $q = T_{qp}(G_{\psi},K)p$.
Suppose now we succeed in tuning $K^\sharp\in \mathscr K$ such that
$(G_{\psi},K^\sharp)$ is $L_2$-stable and satisfies the estimate
$\|T_{qp}(G_{\psi},K^\sharp)\|_\infty < r^{-1}$.  Then by the small-gain theorem
the non-linear loop $(T_{qp}(G_{\psi},K^\sharp),\psi)$ is $L_2$-stable, hence so is
$(G_{nl},K^\sharp)$. This is addressed by the structured mixed $H_\infty/H_\infty$-optimization
program
\begin{eqnarray}
\label{first}
\begin{array}{ll}
\mbox{minimize} & \|T_{qp}(G_{\psi},K)\|_\infty \\
\mbox{subject to}& \|T_{wz}(G,K)\|_\infty \leq (1+\tau) \gamma_\infty\\
&\mbox{$K$ stabilizes $G,G_{\psi}$} \\
&\mbox{$K\in \mathscr K$}
\end{array}
\end{eqnarray}
which optimizes 
stability of the non-linear system $G_{nl}$ under a constraint allowing a controlled loss of performance in the linearized channel $w\to z$, where
$K\in\mathscr K$ ranges over a class of structured controllers in the sense of \cite{AN06a}. 
This is also known as multi-disk optimization \cite{an_disk:05}. 
The algorithmic solution proposed in that reference is implemented in the
{\tt systune} package of \cite{CST2020b}, which we use to solve (\ref{first}) algorithmically.

\begin{proposition}
Suppose the solution $K^\sharp\in \mathscr K$ of {\rm (\ref{first})} satisfies
$\|T_{qp}(G_\psi,K^\sharp)\|_\infty < r^{-1}$. Then the loop $(G_{\rm nl},K^\sharp)$
is stable in the $L_2$-sense. That is, for every $w\in L_2[0,\infty)$
and every $x_0$ the solution of the non-homogenous Cauchy problem
$\dot{x}_{cl} = A_{cl}(K^\sharp)x_{cl} + B_p\phi(C_q x) + B_ww$, $x_{cl}(0)= x_0$ is in $L_2[0,\infty)$.
Moreover $\gamma^\sharp_\infty = \|T_{zw}(G,K^\sharp)\|_\infty\leq (1+\tau)\gamma_\infty$.
\hfill $\square$
\end{proposition}

\begin{remark}
Note that $(G,K)$ and $(G_{\psi},K)$ have different closed loop system matrices.
The $A$-matrix of $G$ is $A$, that of $G_\psi$ is $A_\psi$,  so we have a structured simultaneous stabilization problem,
which is known to be NP-hard for most structures.
\end{remark}

\subsection{Asymptotic sector constraint}
\label{asymptotic}
\begin{figure}[ht!]
\includegraphics[scale=0.8]{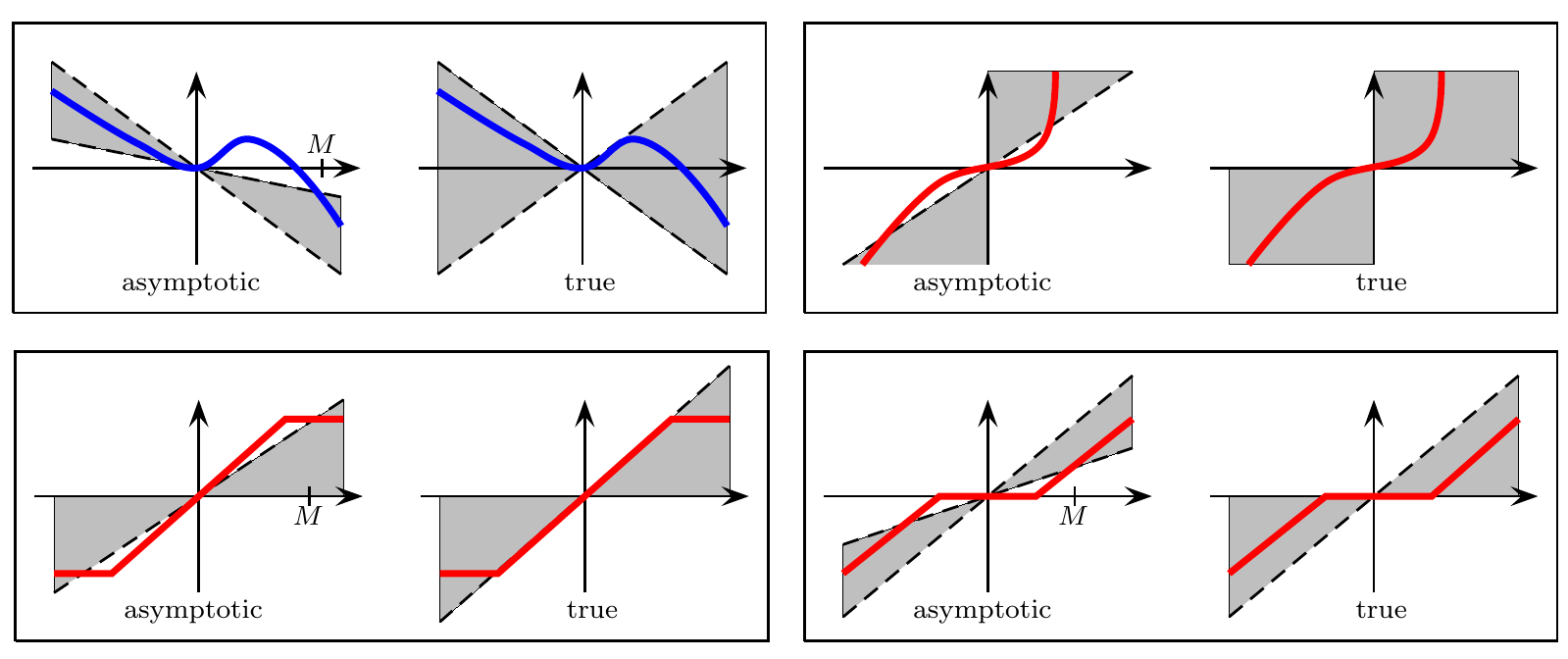}
\caption{Schematic view of true and asymptotic sectors for general sector non-linearity (upper left), positivity (upper right),
saturation (lower left), and dead time (lower right). 
Asymptotic sector constraints are satisfied for large values $|q| >M$.
\label{fig_asymptotic}}
\end{figure}

Apart from the fact that optimization in (\ref{first}) is over structured controllers $K\in \mathscr K$,
the method so far is standard. The situation changes if the sector {\bf sect}$(a,b)$ is  too large, 
so that tuning $K$ to achieve
$\|T_{qp}(G_{\psi},K)\|_\infty  < r^{-1}$ fails. 
Then we have to change strategy!
What we propose in this work is to
choose
a different sector, also noted {\bf sect}$({a},{b})$ for simplicity, which the non-linearity $\phi$ now satisfies only asymptotically, where Fig. \ref{fig_asymptotic} shows schematically what we have in mind.

\begin{definition}
{\bf (Asymptotic sector)}.
{\rm
A non-linearity $p=\phi(t,q)$ satisfies a sector constraint
asymptotically, noted $\phi \sim {\bf sect}(a,b)$, if there 
exist $M,L > 0$ such that for every $t\geq 0$, $(\phi(t,q)-aq)(\phi(t,q)-bq)\leq 0$ for  $|q| > M$, and
$|\phi(t,q)| \leq L$ for $|q|\leq M$. 
}
\end{definition}

Suppose we have identified a new typically smaller sector %
with  $\phi \sim {\bf sect}(a,b)$. (See for instance Fig. \ref{fig_asymptotic}
for some basic examples of asymptotic sector constraints).
We center the non-linearity,
now with the new ${c}=(b+a)/2, {r} = (b-a)/2$, which leads to a new
${\psi}(t,p) = \phi(t,p) - {c}p$, now satisfying
the sector constraint 
$\psi \sim {\bf sect}(-r,r)$ asymptotically.
With $G_\psi$ taken with regard to the new $\psi$, the non-linearity
(\ref{NL1}) is still equivalent to this modified loop
$(G_\psi,\psi)$.
Now
we  consider the mixed {\tt peak-gain}/$H_\infty$-optimization program
\begin{eqnarray}
\label{second}
\begin{array}{ll}
\mbox{minimize} & \|T_{pq}(G_{{\psi}},K)\|_{\rm pk\_gn} \\
\mbox{subject to}& \|T_{wz}(G,K)\|_\infty \leq (1+\tau)\gamma_\infty\\
&\mbox{$K$ stabilizes $G,G_{{\psi}}$} \\
&\mbox{$K\in \mathscr K$}
\end{array}
\end{eqnarray}
where $T_{qp}(G_{{\psi}},K)$ is the channel $p\to q$ of the modified $G_{{\psi}}$ in feedback with $K$.
This optimizes the peak-gain norm of $p\to q$ subject to a controlled loss of $H_\infty$-performance
in the channel $w\to z$ over the optimistic performance 
$\gamma_\infty$ achieved by $K^*\in \mathscr K$.
The algorithmic solution of this novel mixed synthesis  program
will be discussed in Section \ref{new_section}.

We now have the following consequence of the Small-Gain theorem (compare \cite{mareels_hill:92,teel}),
see also Theorem \ref{small_gain}:

\begin{theorem}
\label{theorem1}
Let $K^\flat\in \mathscr K$ be a solution of
program {\rm (\ref{second})} satisfying $\|T_{pq}(G_{{\psi}},K^\flat)\|_{\rm pk\_gn} < r^{-1}$. Then for every input $w\in L_\infty[0,\infty)$ and
all initial conditions $x_0$ the non-linear closed loop
$\dot{x}_{cl} = A_{cl}(K^\flat) x_{cl} + B_p\phi(C_qx) + B_ww$ has trajectories in
$(L_\infty[0,\infty),|\cdot|_\infty)$, and is locally exponentially stable.
\hfill $\square$
\end{theorem}

\newpage

\subsection{Algorithm}
\label{algorithm}
The findings of the previous sections lead to the following strategy:
 \begin{algorithm}[!ht]
 \begin{algorithmic}[1]
\caption{\!\!{\bf :} Mixed peak-gain/$H_\infty$-control of Lur'e system $G_{nl}$}

\State{\bf Steady-state.} 
Compute steady state of non-linear system $G_{\rm nl}$, shift it to origin, and obtain linearization $G$.
\State{\bf Nominal synthesis.} 
Fix performance and robustness specifications and
perform nominal synthesis   for $G$, interpreting non-linearity as a disturbance. Optimal $K^*$ 
gives lower bound $\gamma_\infty = \|T_{wz}(G,K^*)\|_\infty$. 
\State{\bf Sector.}
Using $\phi\in {\bf sect}(a_0,b_0)$, let $c_0=(b_0+a_0)/2$, $r_0=(b_0-a_0)/2$, 
form $\psi_0 = \phi-c_0I$, %
and 
represent non-linear system $G_{\rm nl}$ as loop $(G_{\psi_0},\psi_0)$.
\State{\bf Complementary sector.}
Attempt closed-loop $L_2$-stability by solving
\begin{eqnarray*}
\begin{array}{ll}
\mbox{minimize} & \|T_{qp}(G_{\psi_0},K)\|_\infty \\
\mbox{subject to} & \|T_{wz}(G,K)\|_\infty \leq (1+\tau) \gamma_\infty\\
& K \in \mathscr K
\end{array}
\end{eqnarray*}
If optimal solution $K^\sharp\in \mathscr K$  satisfies $\|T_{pq}(G_{\psi_0},K^\sharp)\|_\infty < r_0^{-1}$, 
quit successfully. Otherwise continue with step 5. 
\State{\bf Asymptotic sector.}
Find asymptotic sector $\phi \sim {\bf sect}(a,b)$. Form ${c}=(b+a)/2$,
${r}=(b-a)/2$, and $\psi = \phi - cI$.
Represent $G_{\rm nl}$ as loop  $(G_{{\psi}},{\psi})$.
\State{\bf Complementary asymptotic sector.}
Attempt closed-loop BIBO-stability in tandem with local exponential stability by solving
\begin{eqnarray*}
\begin{array}{ll}
\mbox{minimize} & \|T_{qp}(G_{{\psi}},K)\|_{\rm pk\_gn}\\
\mbox{subject to} & \|T_{wz}(G,K)\|_\infty \leq (1+\tau) \gamma_\infty\\
&K \in \mathscr K
\end{array}
\end{eqnarray*}
If optimal solution $K^\flat$ satisfies $\|T_{qp}(G_{{\psi}},K^\flat)\|_{\rm pk\_gn} < r^{-1}$ quit successfully.
\end{algorithmic}
\end{algorithm}
 
\begin{remark}
By (\ref{chain}) 
we have $\|G\|_\infty \leq  \|G\|_{\rm pk\_gn}$ even for infinite dimensional systems, so that
$\|T_{pq}(G_{{\psi}},K)\|_{\rm pk\_gn} < 1$ implies
$\|T_{pq}(G_{{\psi}},K)\|_\infty < 1$. Therefore 
it makes no sense to choose the asymptotic sector as a true sector.
We need $\phi \sim {\bf sect}(a,b)$, but must have
$\phi \not\in {\bf sect}(a,b)$, as $\phi\in {\bf sect}(a,b)$   would mean trying step 4  again, 
saddled with the even harder
constraint $\|T_{pq}(G_{{\psi}},K)\|_{\rm pk\_gn} < 1$.
\end{remark}

\subsection{Best asymptotic sector}
\label{best}
Working with asymptotic sectors offers additional flexibility over conventional sectors, 
which we now exploit.
Consider
step 6 of the algorithm.   
Instead of choosing the asymptotic sector {\bf sect}$(a,b)$,
which is the same as choosing $c,r$,
we could in the first place only choose ${c}$. With $\psi(t,q) = \phi(t,q) - cq$,
we solve program %
\begin{align*}
&(\ref{second}') \hspace{4cm}
\begin{array}{ll}
\mbox{minimize} & \|T_{pq}(G_{\phi-c},K)\|_{\rm pk\_gn} \\
 \mbox{subject to}& \|T_{wz}(G,K)\|_\infty \leq (1+\tau) \gamma_\infty\\
&K\in \mathscr K
\end{array}
\hspace{4cm}
\end{align*}
where $r$ is not yet determined.
The optimal controller $K({c})\in \mathscr K$ now depends on ${c}$, and as in step 6 of the algorithm, 
provides the value
\begin{equation}
\label{1/r}
{r}({c}) := 1/ \|T_{pq}(G_{\phi-{c}},K({c}))\|_{\rm pk\_gn}.
\end{equation}
This gives
a curve ${r}={r}({c})$, and on putting
$a=c-r(c)$, $b=c+r(c)$, it remains to check whether
{\bf sect}$(a,b)$ is an asymptotic sector for $\phi$. 

An interesting case is when $\phi$
has a slope at infinity, i.e., when
$\lim_{|x|\to\infty} \frac{\phi(t,x)}{x} = q_\infty$ exists independently of $t$.
Then every choice $a < q_\infty < b$ gives $\phi \sim {\bf sect}(a,b)$,
so with (\ref{1/r}) we have to check whether
\begin{equation}
\label{works}
{c} - {r}({c}) < q_\infty < c+r(c).
\end{equation}
As soon as this holds, 
we have {\it a posteriori} found an asymptotic sector $\phi \sim {\bf sect}(a,b)$  as
$a = {c}-{r}({c})$ and $b={c}+{r}({c})$. 
We can then also  determine
the parameters $L,M$ in the definition of an asymptotic sector. 
The smaller $M$, the closer the asymptotic sector comes to a true sector.

Since the asymptotic sector is chosen in response to
failure of the true sector {\bf sect}$(a_0,b_0)$, we typically initialize the search for $a,b$ 
by values ${c}$ close to ${q}_\infty$, as this increases the chances of program ($\ref{second}'$)
to succeed.
This highlights why
an asymptotic sector typically will not satisfy $0\in [a,b]$, whereas this is always satisfied for the true sector. 

\begin{remark}
For every $c$  the aperture $2r(c)$ of the candidate sector is maximized through program (\ref{second}) due to (\ref{1/r}). Over the range of
those $c$ where (\ref{works}) holds, the resulting curve $(c,r(c))$ serves as a Pareto optimal front, from which we will pick our ultimate $c$.
The decision will  
not just be based on the size of the aperture of {\bf sect}$(a,b)$,
it may also matter how close {\bf sect}$(a,b)$  is to a true sector, how large the constant $k=LM$ is, and ultimately,
how $G_{\rm nl}$ in loop with $K(c)$ 
behaves in non-linear simulations. This method will be applied to the slipstick study in Section \ref{sect_slipstick}.
\end{remark}

\section{Application: BIBO-stable control of slipstick vibrations}
\label{sect_slipstick}
We consider control of a damped wave equation with instability caused by non-linear
boundary anti-damping dynamics, 
\begin{align}
\label{slipstick}
\begin{split}
x_{tt}(\xi,t) &= x_{\xi\xi}(\xi,t) - 2 \lambda x_t(\xi,t) , \;\; 0 < \xi < 1, t\geq 0 \\
G_{\rm nl}: \qquad x_\xi(1,t) &= -x_t(1,t) + u(t) \\
\alpha x_{tt}(0,t) &= x_\xi(0,t) + qx_t(0,t) + \phi(x_t(0,t))
\end{split}
\end{align}
where %
$(x,x_t)$ is the state, $u(t)$ the boundary control, and the measured outputs are
\begin{equation}
\label{outputs}
y_1(t) = x_t(0,t), \quad y_2(t) = x_t(1,t).
\end{equation}
The non-linearity satisfies $\phi(0)=0$, $\phi'(0)=0$, so that the linearized
system $G$ in (\ref{L1}) is obtained by dropping the term $\phi(x_t)$.
System $G_{\rm nl}$ has among others been used to model slipstick vibrations in drilling systems, see \cite{an1} and the references given there. The challenge is to
control $G_{\rm nl}$ with a finite-dimensional controller $u=Ky$  of simple, implementable structure
such that slipstick caused by the non-linear boundary friction term $\phi(x_t)$ can be avoided or at least mitigated. 

In these applications
$\lambda \geq 0$, $\alpha \geq 0$,  $q\geq 0$ are typically positive, and
the non-linearity derives from a frictional force depending on the angular velocity of the drill
$$
T(\omega) = \gamma_1 \omega + \left(\gamma_2 + \gamma_3 e^{-\gamma_4|\omega|} \right){\rm sign}(\omega),
$$
exhibiting a sharp jump at $\omega=0$, which
based on experimental evidence in comparable situations \cite{ravanbod},
is slightly mollified around $0$. 
With $\bar{\omega} >0$, 
the nominal angular speed of the drill, step 1 of the algorithm leads to the centered non-linearity
\[
\phi(\omega) = T(\bar{\omega}) - T(\bar{\omega}+\omega) + T'(\bar{\omega})\cdot \omega
\]
as in (\ref{NL1}), 
shown in
Fig. \ref{sectors} for two of the scenarios studied in \cite{an1}.
This corresponds to step 1 of the algorithm.

\vspace{.2cm}
\begin{center}
\includegraphics[scale = 0.42]{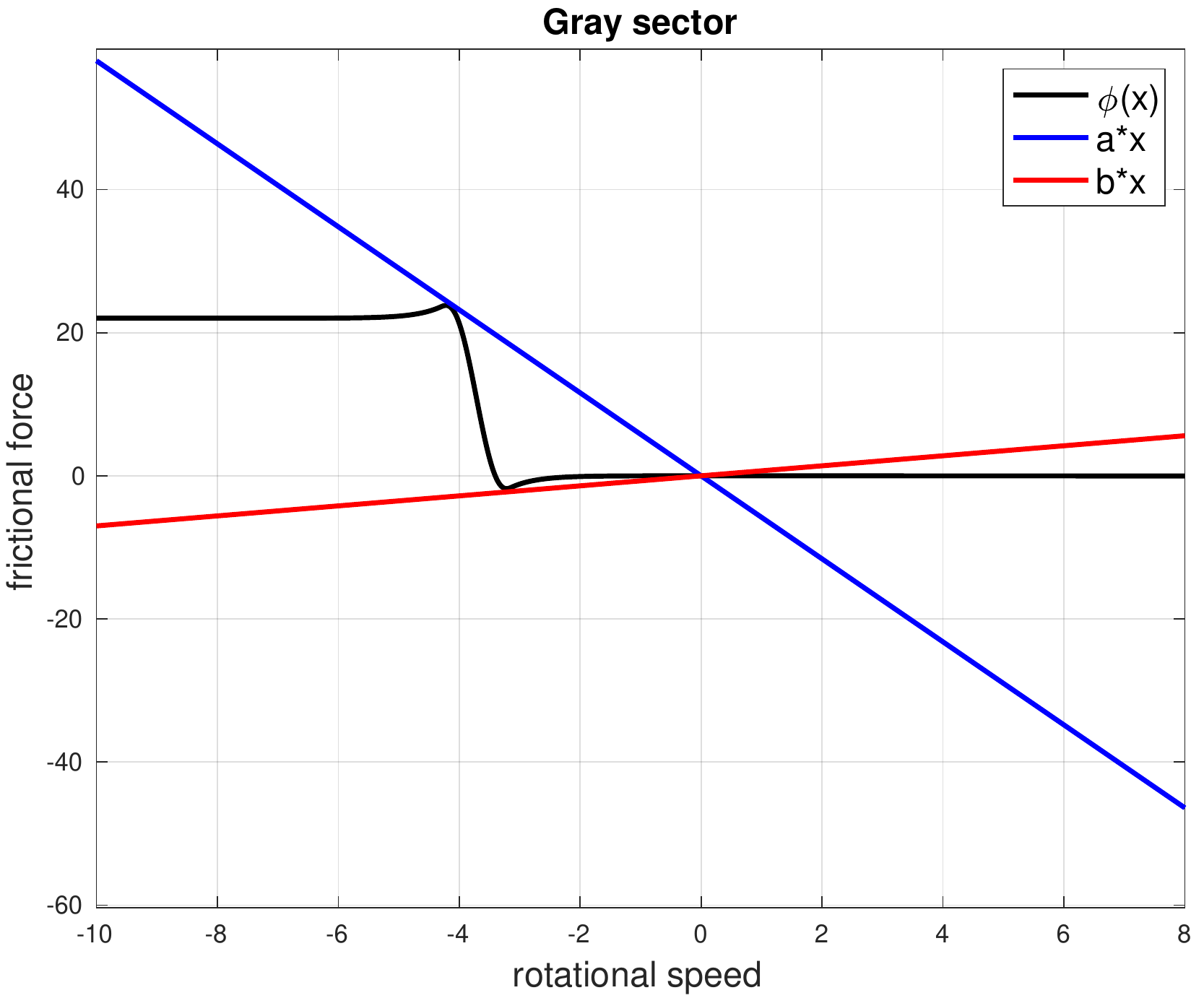} $\quad$
\includegraphics[scale = 0.42]{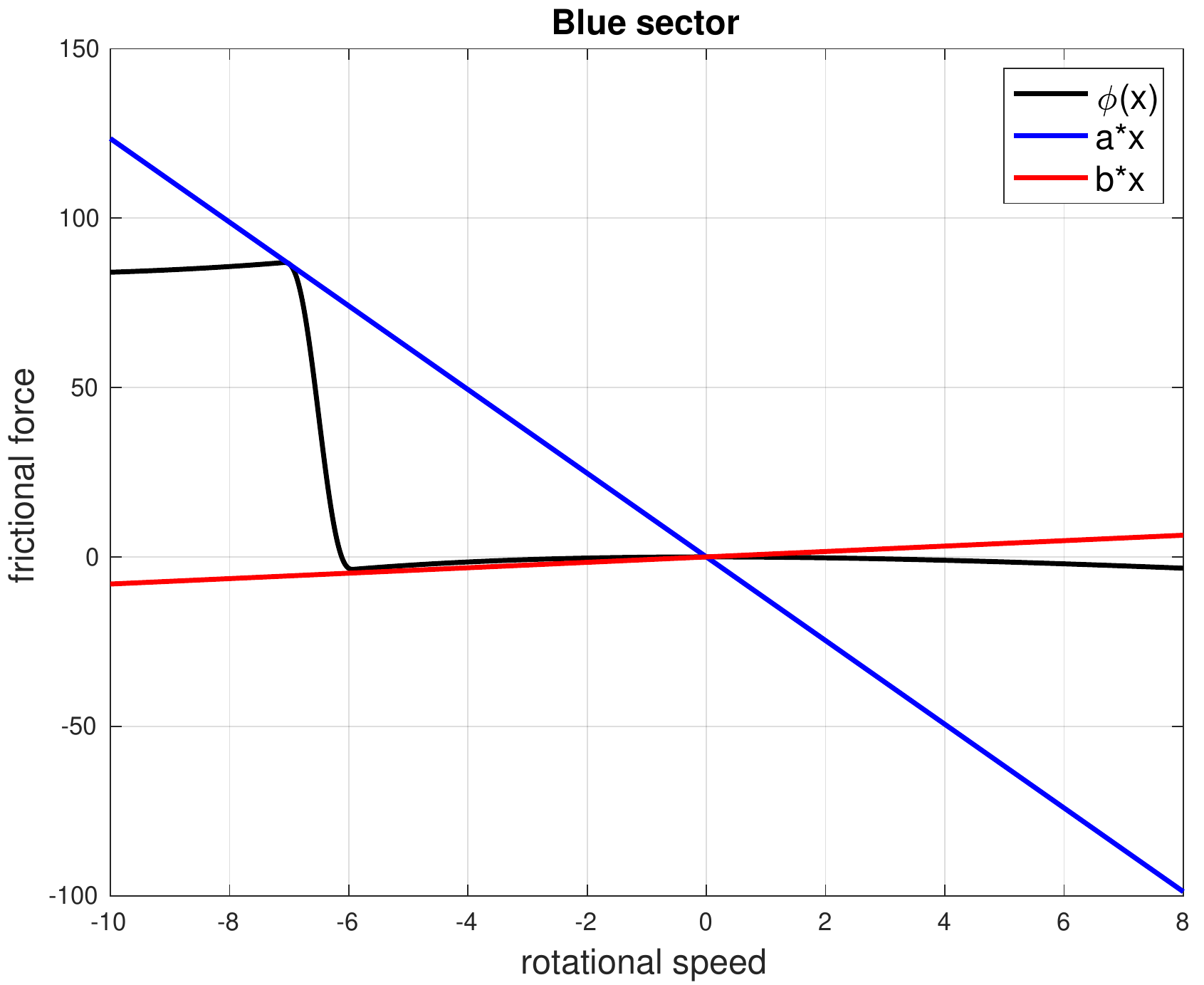}
\end{center}
\vspace*{-1cm}
\begin{figure}[ht!]
\caption{Sector non-linearity $\phi$ in slipstick model with strong jump
at $-\overline{\omega}$ for two scenarios labeled 'gray' and 'blue'  in \cite{an1}.  \label{sectors}}
\end{figure}

We continue to follow the pattern of the algorithm. 
By \cite[sect. 3-4]{an1}
the number of unstable open-loop poles of $G$ is $n_p(q,\alpha,\lambda) \in \{0,1,2\}$,
and several scenarios 'gray', 'blue', 'red', 'magenta' and 'green' were analyzed.  The blue scenario, on which we focus here,
has two unstable open-loop poles for the chosen $\lambda$, the numerical parameters gathered in Table \ref{table1}.  Here the aperture of the sector 
{\bf sect}$(a_0,b_0)$ in Fig. \ref{sectors} (right) is extremely large, and
step 4 of the algorithm  fails even when a rather conservative $\tau$ is chosen. 
This is where we use the  asymptotic sector
of step 5 of the algorithm.

\vspace{.2cm}
\begin{table}[ht!]
\title{Numerical values for slipstick study}
\begin{center}
\begin{tabular}{|| c || c | c | c | c | c | c | c | c   ||}
\hline\hline
 & $q$ & $\alpha$ & $\lambda$ & $\overline{\omega}$  &$\gamma_1$ & $\gamma_2$ & $\gamma_3$ & $\gamma_4$ \\
\hline
gray & 0.0019 & 0.7994 & 0.1957 & 3.7186 &1.002e-4 & 11.0034 &6.6020&2.4203\\
\hline
blue & 0.9797 & 0.1828 & 0.5477 & 6.5044 &1.002e-4 &28.8697 &17.3218 &0.1537\\
\hline \hline
\end{tabular} 
\end{center}
\caption{ \label{table1}}
\end{table}

\begin{figure}[!ht]
\includegraphics[width=0.49\textwidth]{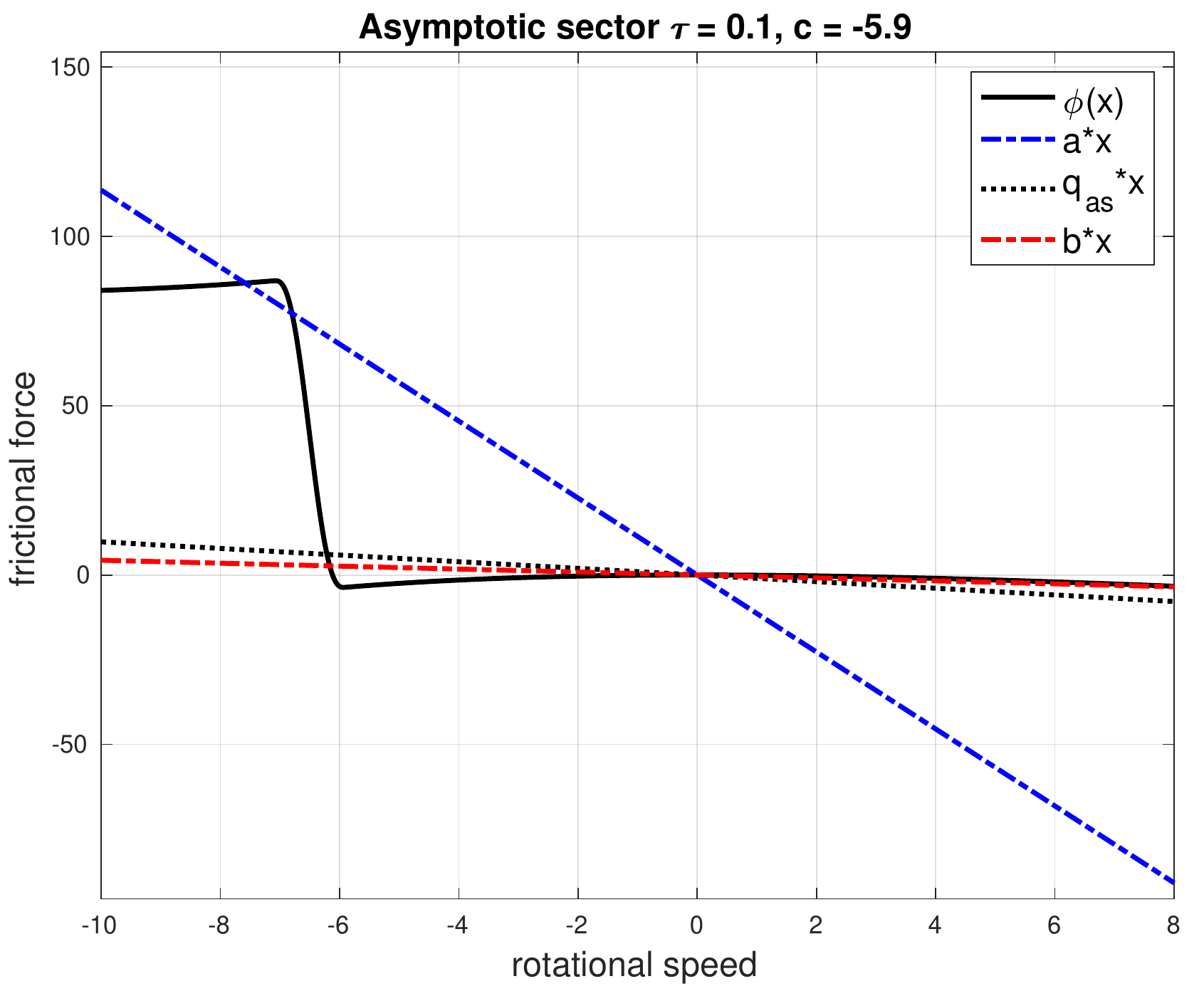}
\includegraphics[width=0.49\textwidth]{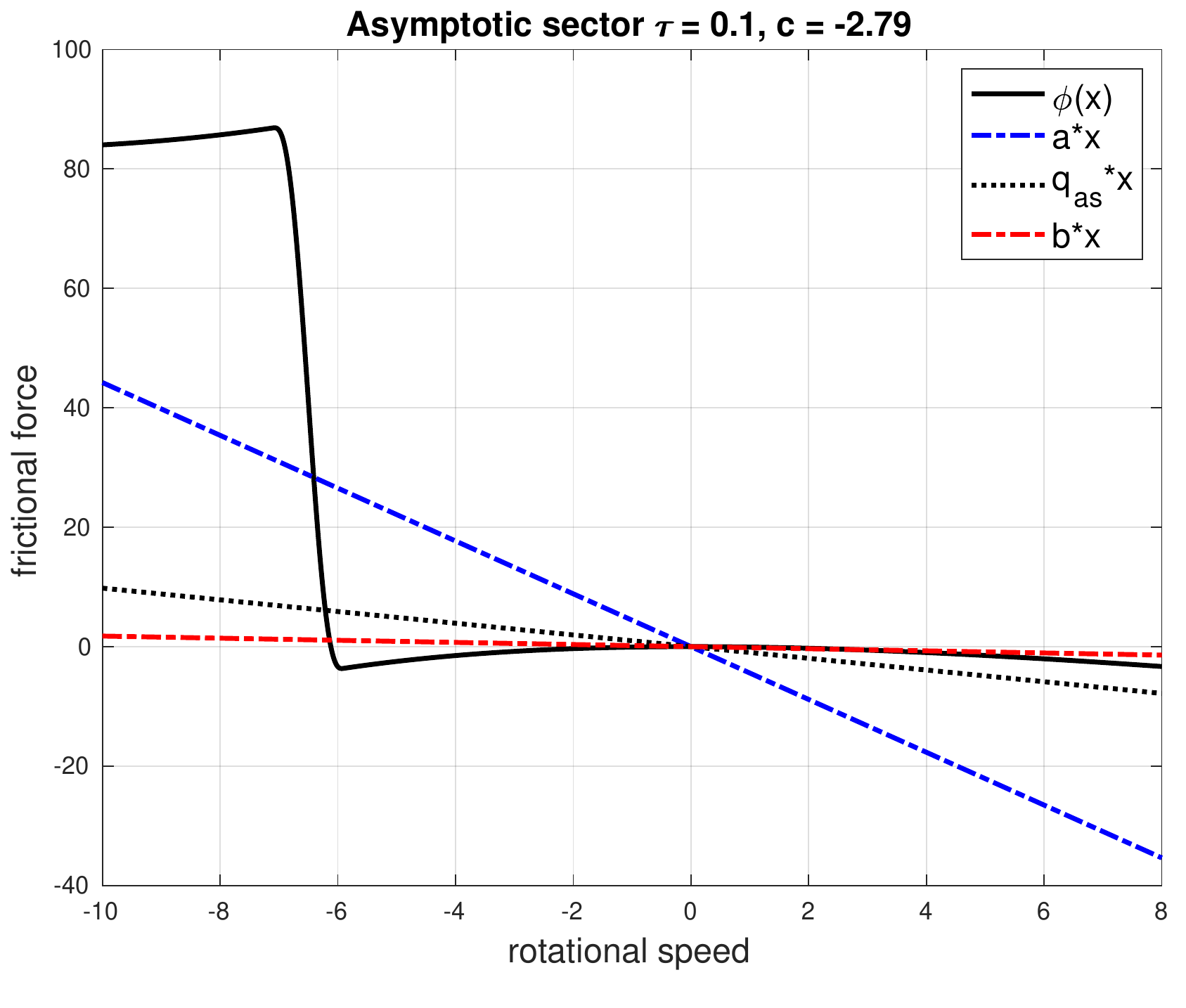}
\caption{Asymptotic sectors for $\tau = 0.1$ and $c=-5.9$, $c=-2.79$, with $r(c)$ computed via (\ref{opt}).
Constraints are satisfied for large angular velocities.\label{as_sectors}}
\end{figure}

From \cite[Lemma 4]{an1} we know that $\phi(\omega)$ behaves asymptotically
as $\phi(\omega) \sim a_\pm + \phi'(\infty) \omega$ for $\omega \to \pm \infty$ with $\phi'(\infty) = -0.9797$. 
This means, every choice $a < \phi'(\infty) < b$ gives rise to an asymptotic sector
$\phi \sim {\bf sect}(a,b)$. This can for instance be seen in Fig. \ref{as_sectors}.
We now have to give the details of step 6 of the algorithm.

As in Fig. \ref{optimistic} (right), we consider a nominal linear model, where the non-linearity
is interpreted as a disturbance $w$. We
optimize the closed-loop $H_\infty$-channel $T_{(r,w) \to (z_u,z_y)}(G,K)$, where a high pass filter $W_u$ is used for the control signal,
and a low-pass for tracking of output $y_1$, which corresponds to the rotational speed at the drill bit. The rationale is that attenuating the disturbance
$w$ should reduce the effect of the non-linearity.
This leads to the nominal $H_\infty$-performance $\gamma_\infty = \|T_{(r,w) \to (z_u,z_y)}(G,K_\infty)\|_\infty$ for an
$H_\infty$-controller $K_\infty \in \mathscr K$. As proved in \cite{an1}, exponential stabilizability
and detectability of the linear open loop guarantee that the linear closed loop $T_{(r,w) \to (z_u,z_y)}(G,K_\infty)$  is not only $H_\infty$-stable, but
even exponentially stable, and as a consequence, BIBO-stable.

\begin{figure}[ht!]
\includegraphics[scale=0.35]{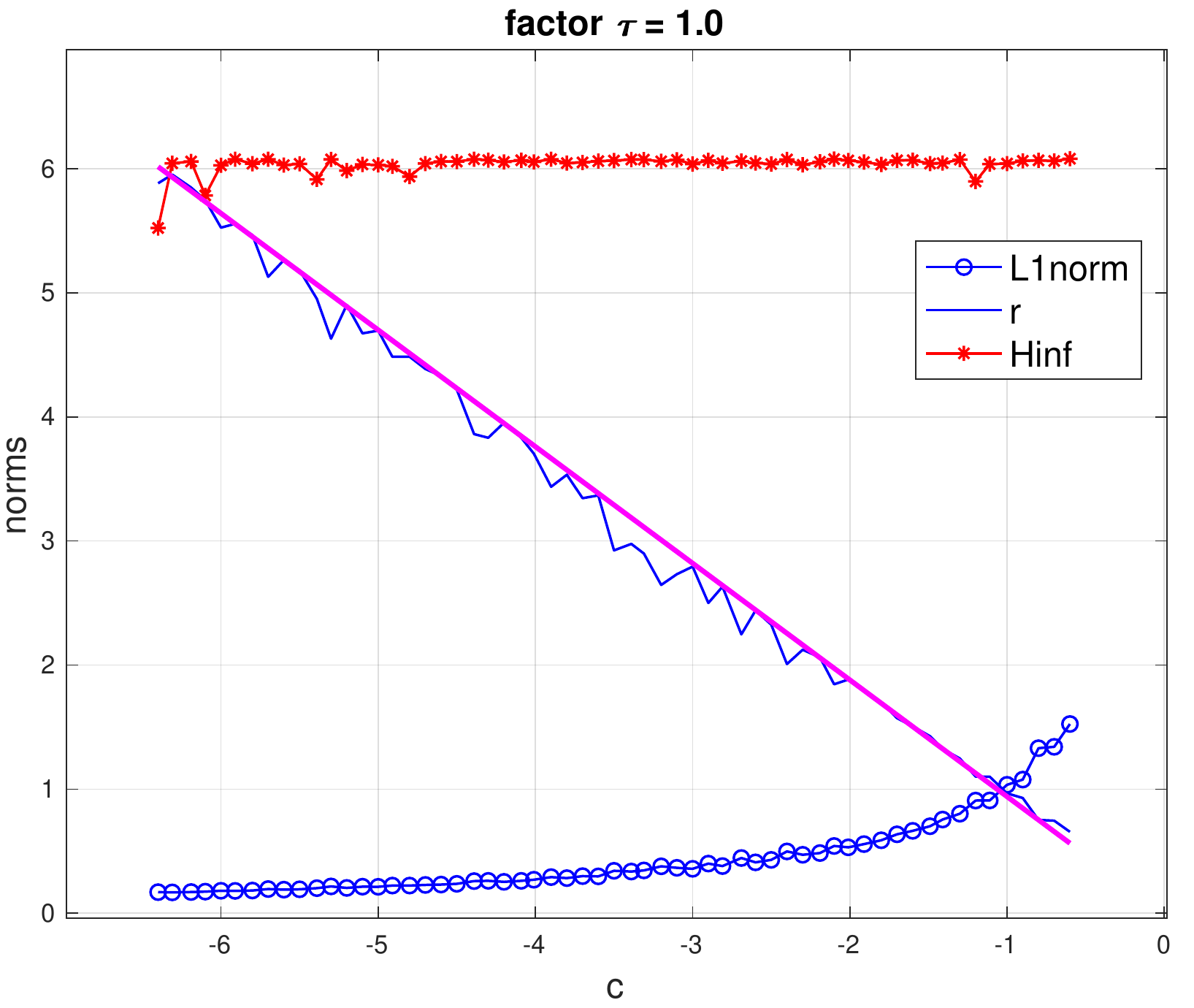}
\includegraphics[scale=0.35]{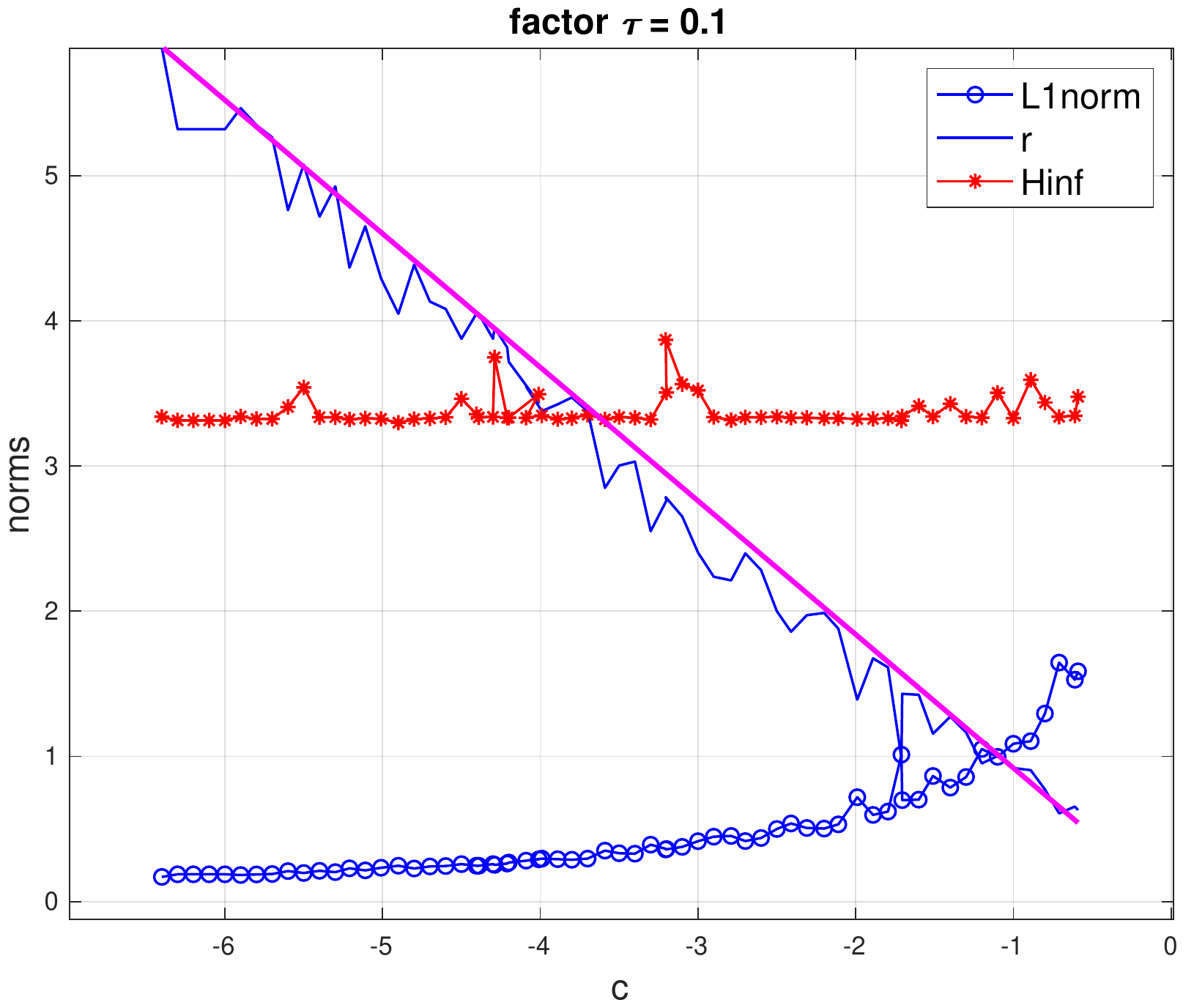}
\includegraphics[scale=0.35]{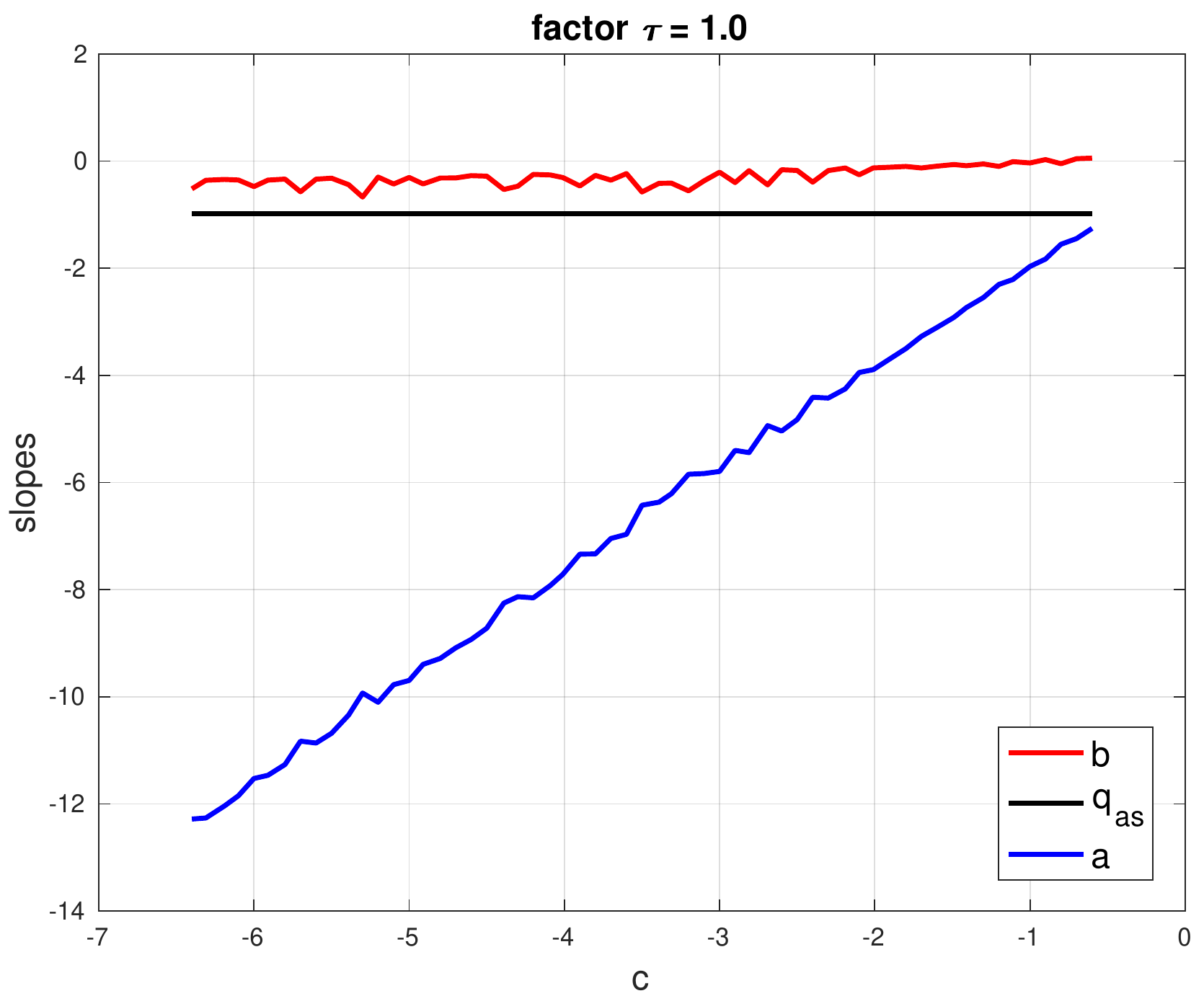}
\includegraphics[scale=0.35]{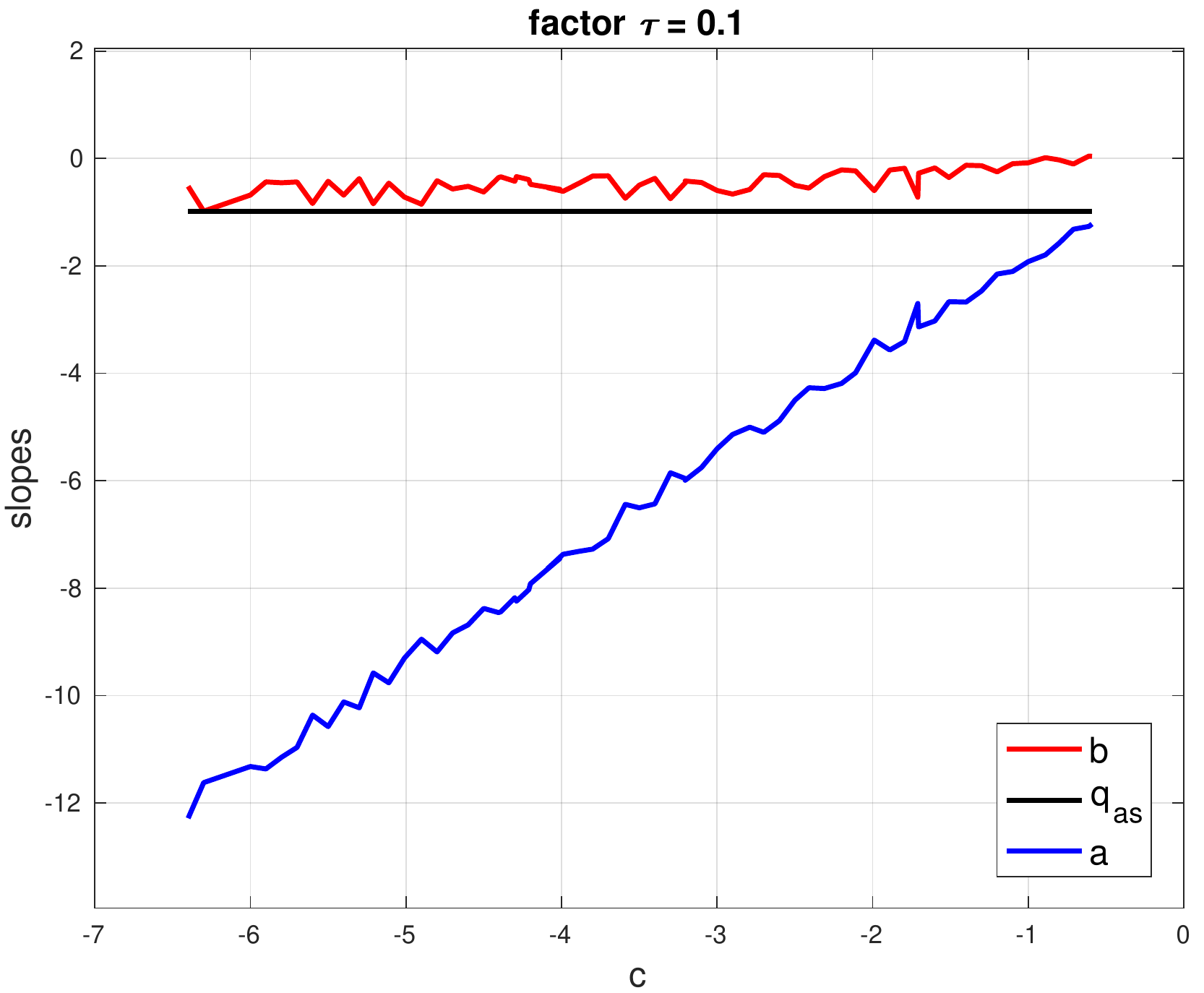}
\caption{Program (\ref{opt}) was run for $c\in [-6.2,-0.5]$ and $\tau = 0.1$ (right),  $\tau= 1.0$ (left). Upper row shows
best achieved $1/r(c)$ as ${\color{blue}-{\rm o}-}$,
with $r(c)\approx -0.94\cdot c$ shown as ${\color{blue}-\!-}$ for $\tau = 1.0$ and $-0.92$ for $\tau = 0.1$, depending essentially linearly on $c$ (magenta).  Nominal  $H_\infty$-norm  
is $3.03$, and ${\color{red}-\ast-}$ shows result after optimization. Lower row shows slopes $a$ (blue), $b$ (red) of asymptotic
sectors with
$a < q_{as} < b$ and $b = c+r(c)$, $a=c-r(c)$. 
\label{fig4}}
\end{figure}

Since it is necessary
to consider asymptotic sectors, 
we choose the parameter $c$, and obtain the loop transformed  system $G_{\phi-c}$. This corresponds to
representing the non-linearity as a feedback loop as in Fig. \ref{optimistic} (left), and we now have to optimize the peak gain norm
of the closed-loop channel
$p \to q$ in $G_{\phi-c}$. This is the mixed peak-gain/$H_\infty$ program
\begin{eqnarray}
\label{opt}
\begin{array}{ll}
\mbox{minimize} &  \|T_{qp}(G_{\phi-c},K)\|_{\rm pk\_gn} \\
\mbox{subject to}&  \|T_{(r,w)\to (z_u,z_y)}(G,K)\|_\infty\leq (1+\tau)\gamma_\infty\\
&\mbox{$K$ stabilizes $G,G_{\phi-c}$}\\
&K\in \mathscr K
\end{array}
\end{eqnarray}
which
leads to the optimal $K_{\rm opt}(c)\in \mathscr K$ with value $1/r(c) =  \|T_{qp}(G_{\phi-c},K_{\rm opt}(c))\|_{\rm pk\_gn}$.
Then we compute
$a = c-r(c)$, $b= c + r(c)$, and if $a < q_{\infty} < b$, then the non-linearity is asymptotically in the sector,
so that the non-linear closed loop is BIBO-stable. 
In our experiment $\mathscr K$ designates 3rd order controllers, which due to $n_y=2$, $n_p=1$ leads to $18$ 
optimization variables.
Fig. \ref{fig4} shows these curves for two scenarios $\tau=0.1$ and $\tau=1.0$, with $c$ in the range 
$c\in [-6.2,-0.5]$.

\begin{figure}[ht!]
\includegraphics[width=0.49\textwidth]{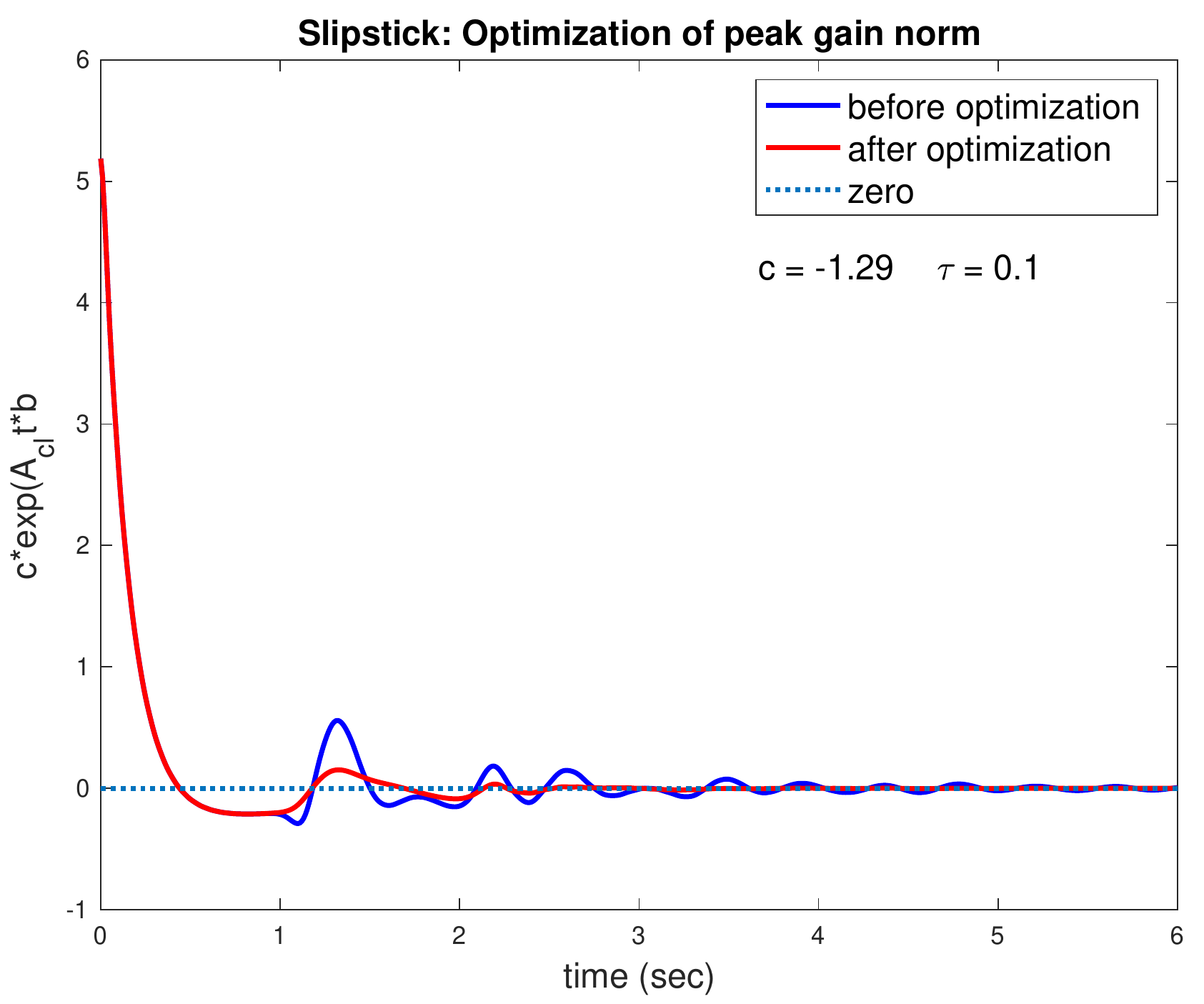}
\includegraphics[width=0.49\textwidth]{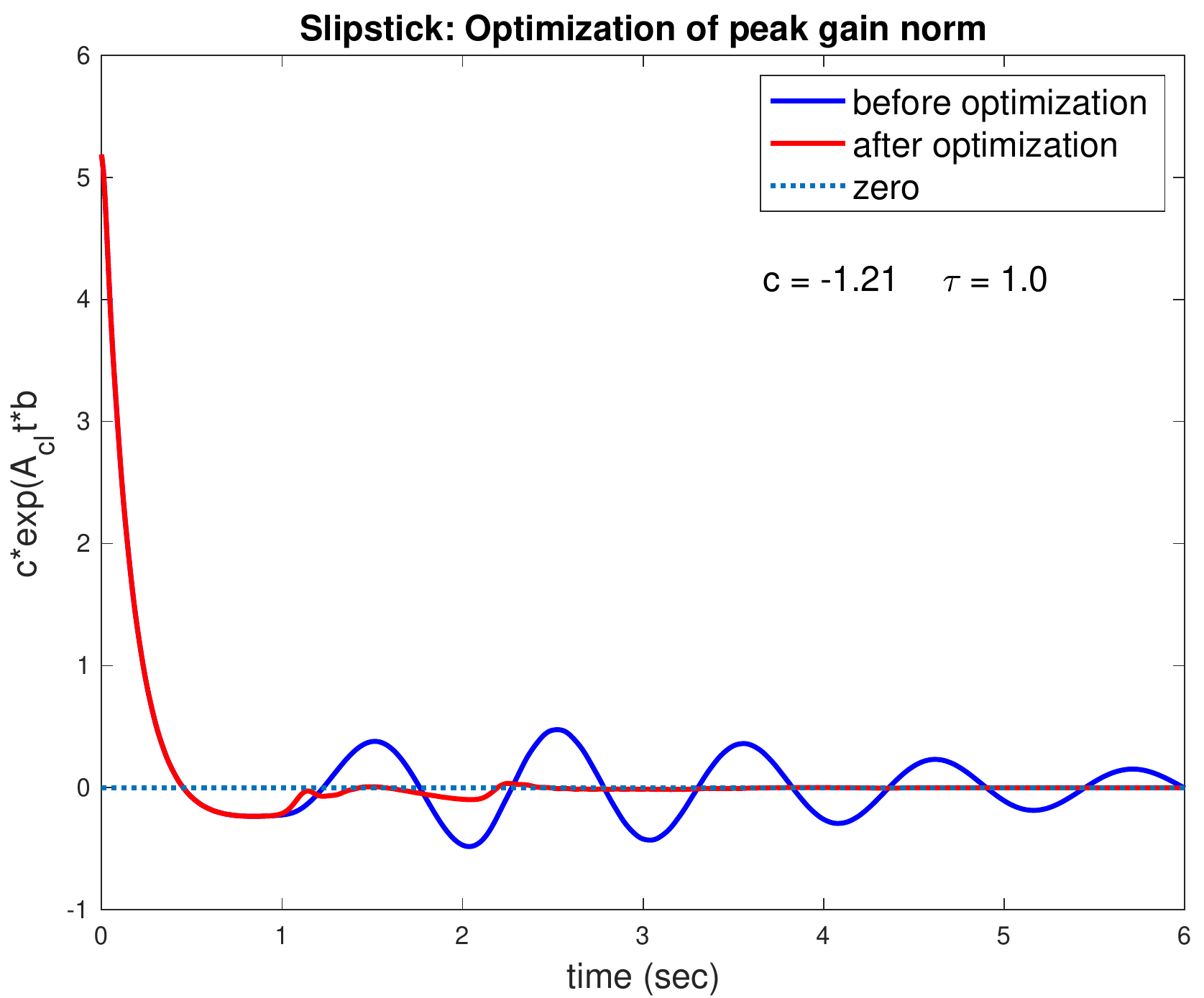}
\caption{Impulse response before and after optimization (\ref{second}) for two cases $c=-1.29$, $\tau = 0.1$ and $c=-1.21$, $\tau=1.0$. \label{responses}}
\end{figure}

\begin{remark}
Application of our theory requires some preparation, as the system is now infinite dimensional
and of boundary control type.  
We have to clarify the meaning of the impulse
response $ce^{A_{cl}t}b$ representing the closed-loop channel $w\to z$. 
Following \cite[Sect. 3.3]{CurtainZwart:1995}, \cite[sect. 5]{an1}, the
linear wave equation and boundary feedback controller can after a change of variables be represented as an abstract boundary control system
\begin{align*}
\dot{x} &= \mathscr Ax, \;
\mathscr Px = u+w\\
y &= \mathscr Cx, \; z = \mathscr C_1 x\\
u &= Ky
\end{align*}
where $u=Ky$ is a finite-dimensional controller,  and 
where $\mathscr A: D(\mathscr A)\to Z$, $Z$ a separable Hilbert space, $u(t)\in \mathbb R^p$, $\mathscr P:D(\mathscr P)\to \mathbb R^p$,
$D(\mathscr A)\subset D(\mathscr P)\subset Z$, $D(A) = D(\mathscr A) \cap {\rm ker}(\mathscr P)$ is dense and
$A = \mathscr A|D(A)$ generates a $C_0$-semi-group on $Z$. Moreover, there exists a bounded operator $B\in L(\mathbb R^p,Z)$ such that
$Bu\in D(\mathscr A)$ for every $u$, $\mathscr A B\in L(\mathbb R^p,Z)$, and $\mathscr PBu=u$ for every $u$.  In \cite[Thm. 2]{an1} the case
$w=0$ was handled, and in order to accommodate the Lur'e non-linearity, we have to make a slight extension. As \cite[Thm. 2]{an1} shows, a
finite-dimensional 
$H_\infty$-stabilizing controller with minimal representation renders the closed loop in this
state-space representation exponentially stable. That means the channel $w \to z$ is represented as
$ce^{A_{cl}t}b$, where $b(\xi) w\in D(A_{cl})$ and $A_{cl}$ generates an exponentially stable semi-group. Hence $e^{A_{cl}t}b$ is a classical
solution, and since $A_{cl}$ is exponentially stable, $e^{A_{cl}t}b \in L^1$ by the Datko-Pazy theorem \cite[Thm. V.1.8]{engel_nagel}. This implies
$ce^{A_{cl}t}b\in L^1$.  In consequence, the impulse response is convenient to optimize, even though we expect a singularity
at $t=0$ (see e.g. Fig \ref{responses}). 
\end{remark}

\begin{figure}[ht!]
\includegraphics[scale=0.32]{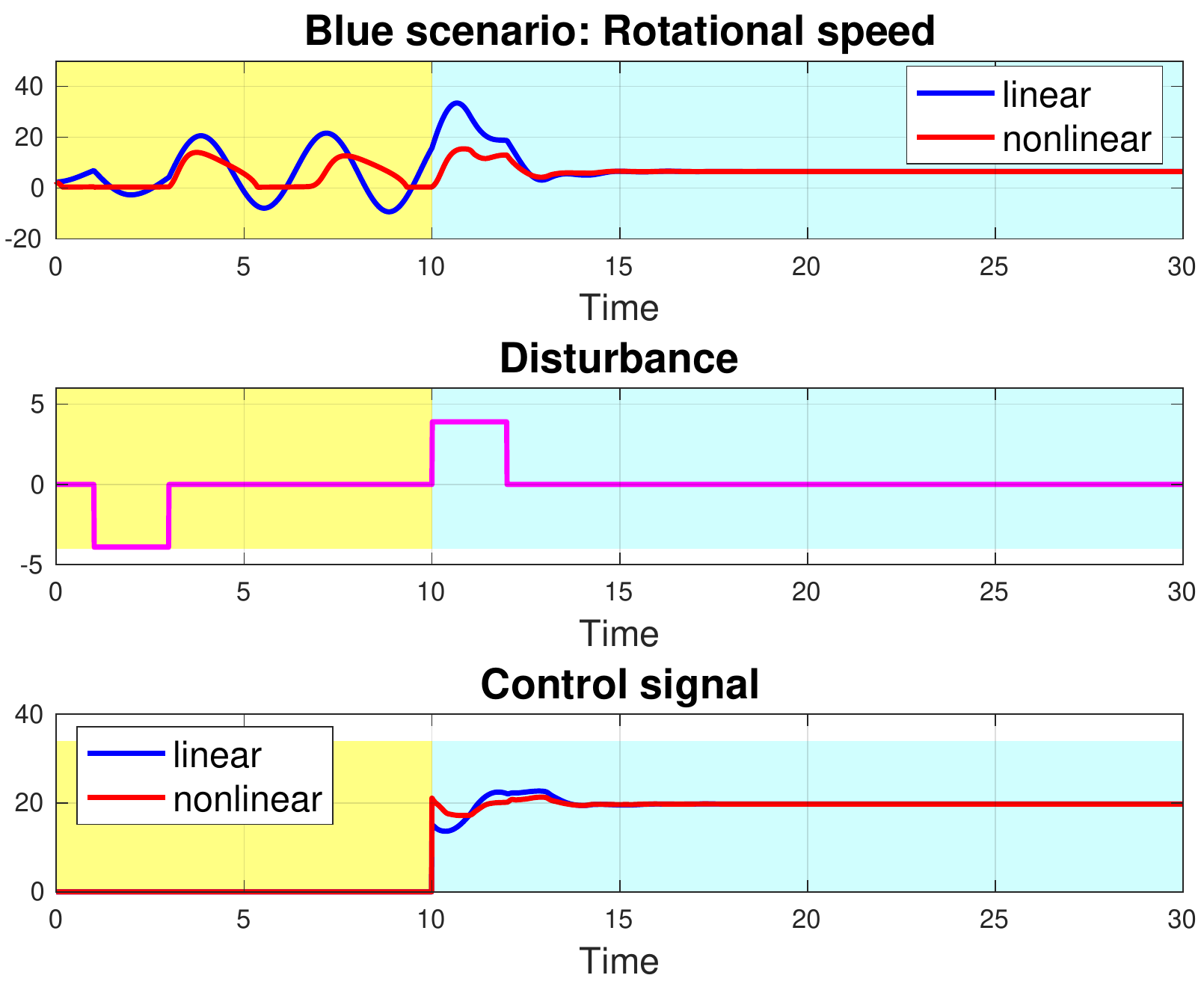}
\includegraphics[scale=0.32]{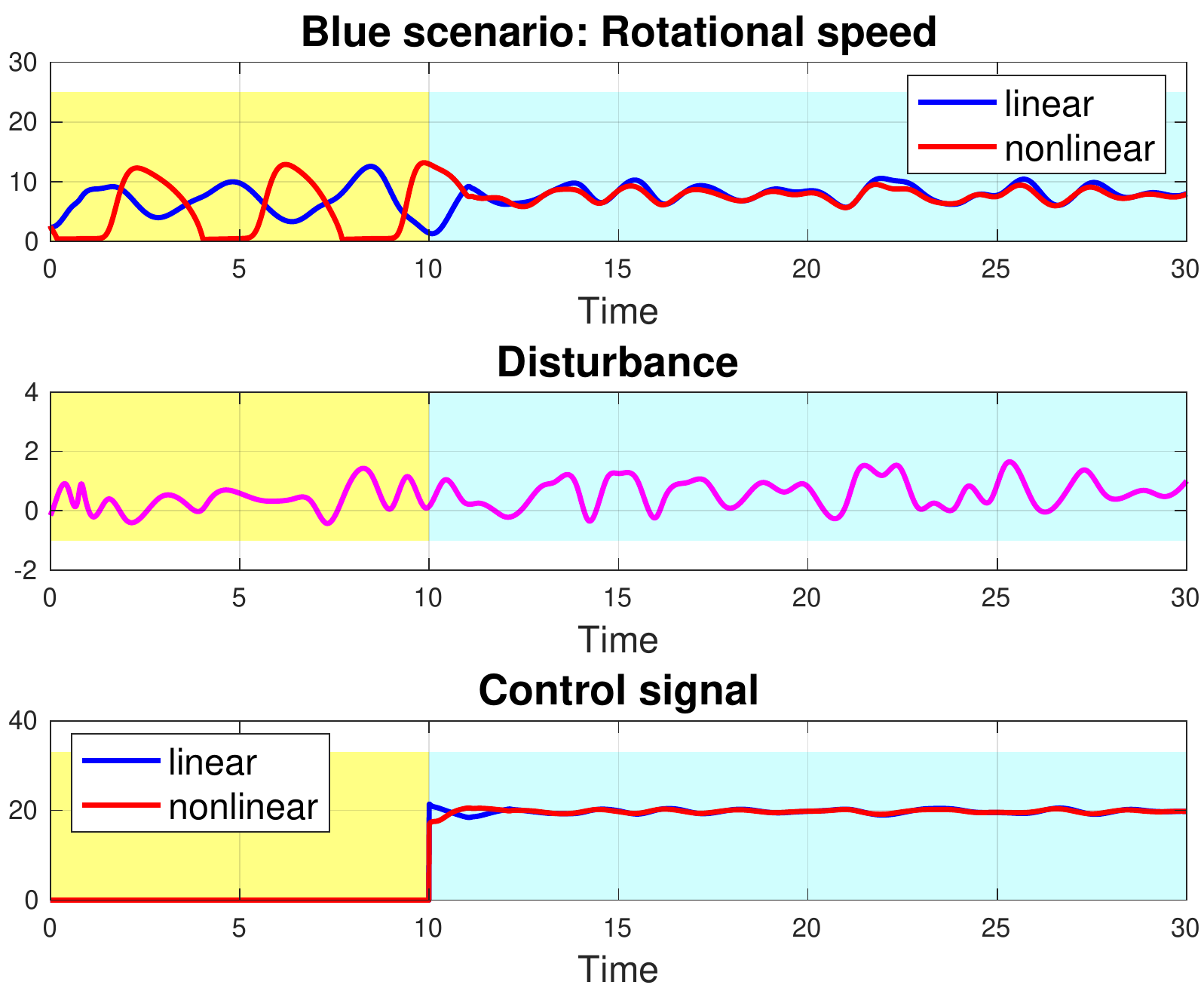}
\includegraphics[scale=0.32]{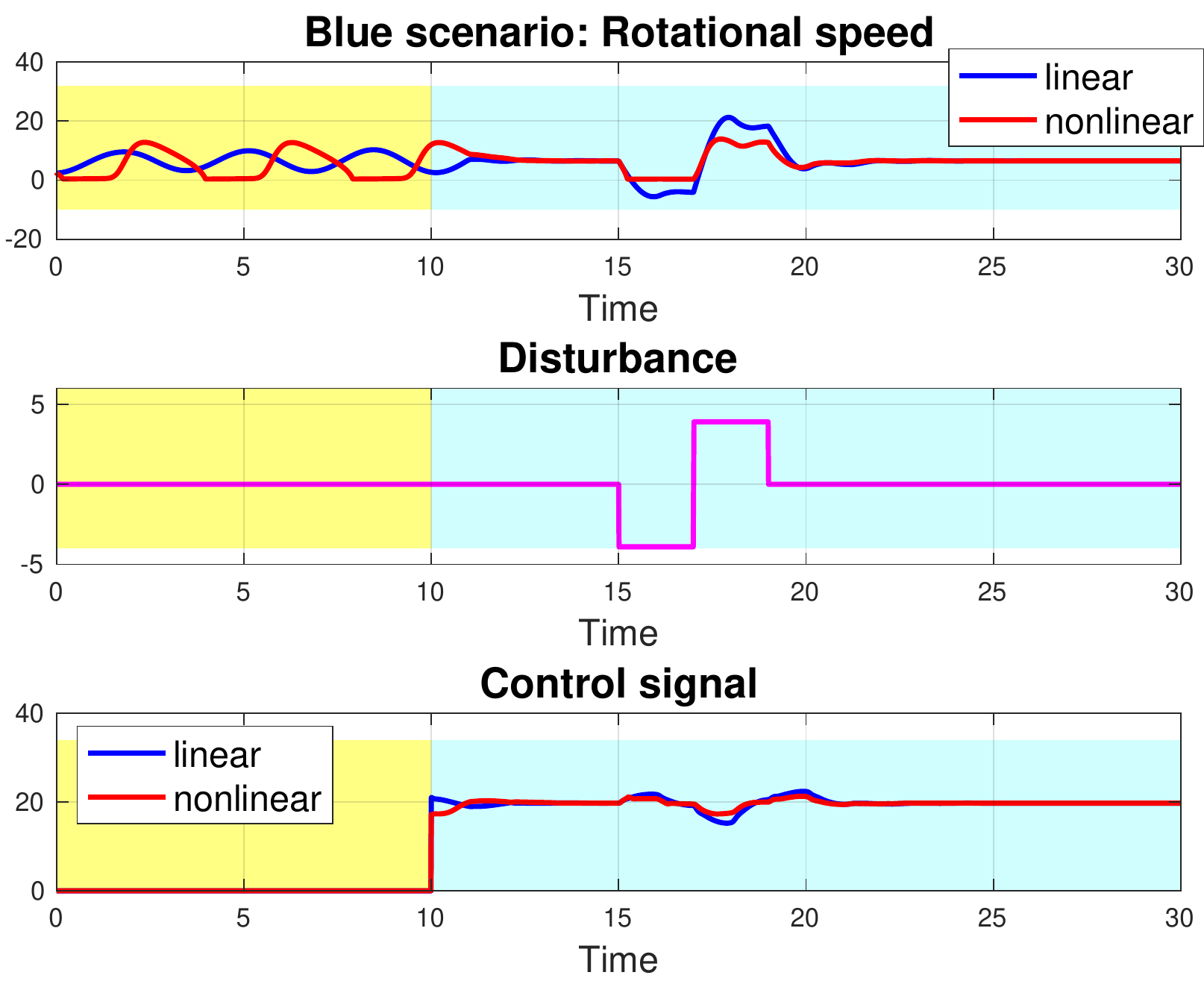}

\caption{$c=-5.21$, $\tau = 0.1$. Initial value below steady state, causing slipstick. Controller switched on at $t=10$. Uncontrolled system shows slipstick. \label{fig7}}
\end{figure}

\begin{figure}[ht!]
\includegraphics[scale=0.318]{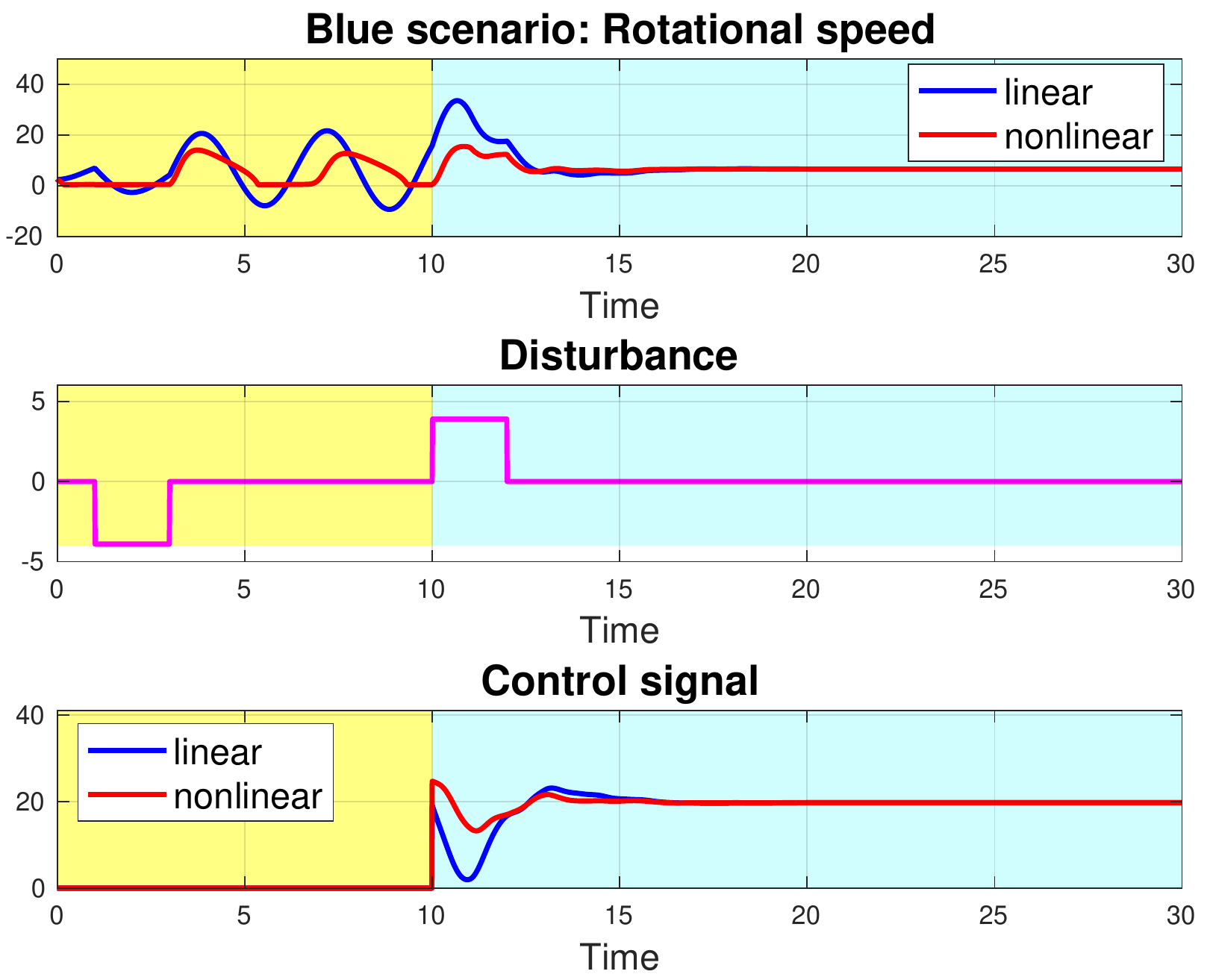}
\includegraphics[scale=0.318]{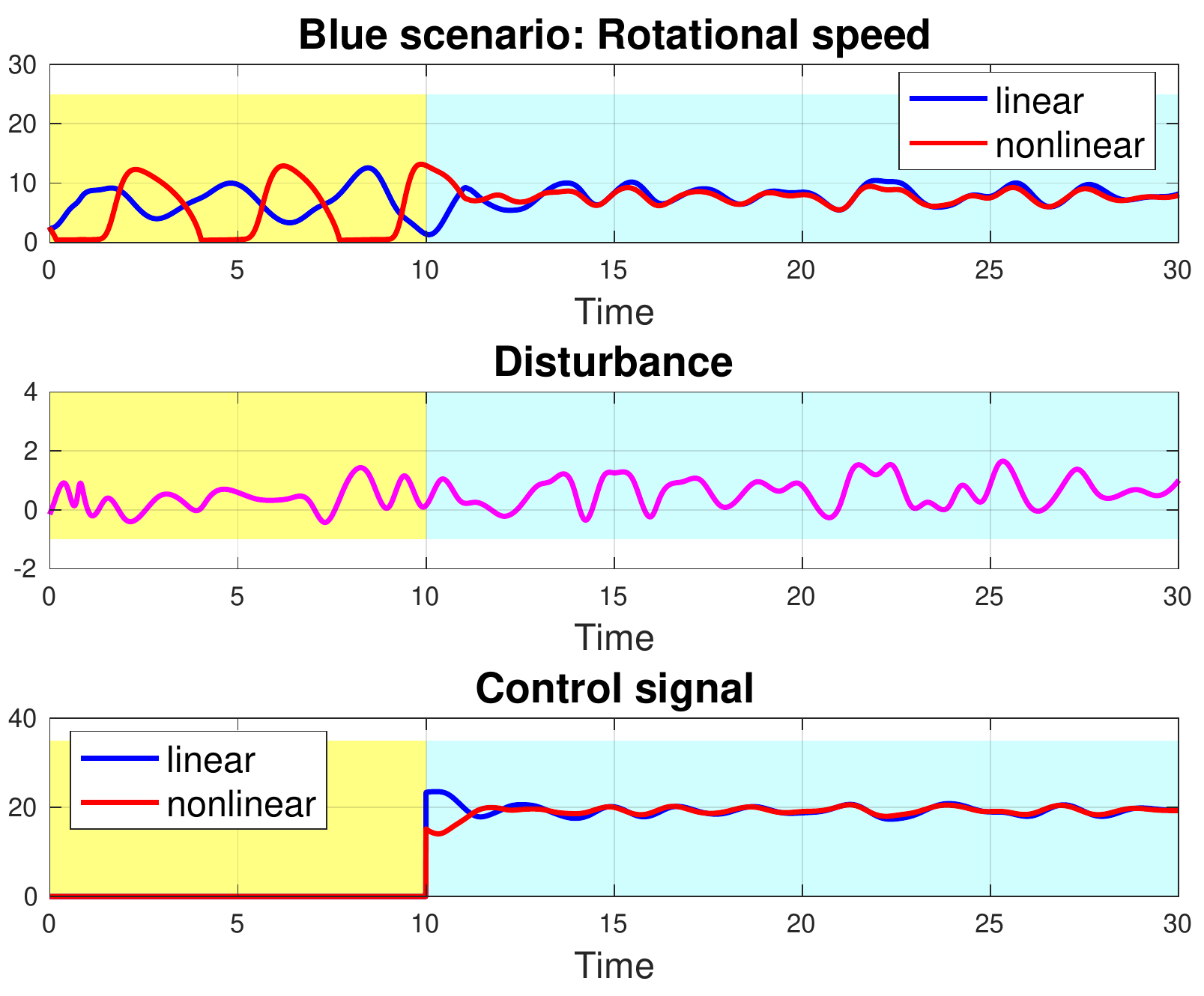}
\includegraphics[scale=0.318]{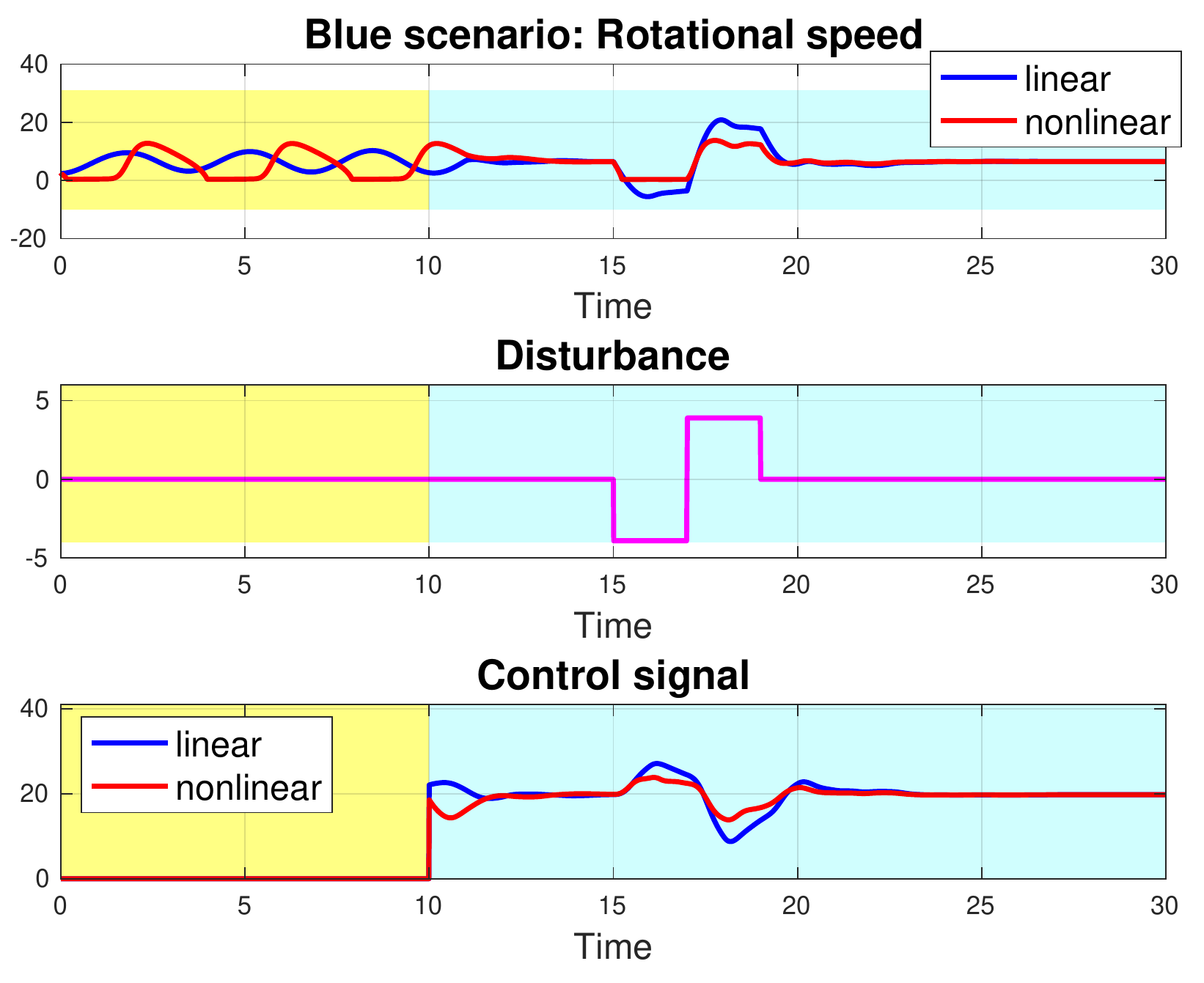}

\caption{$c=-2.79$, $\tau = 0.1$. Initial value below steady state, causing slipstick. Controller switched on at $t=10$. Three types of disturbances. Controller manages to free system from slipstick. \label{fig8}}
\end{figure}

\begin{remark}
We recall that for linear systems BIBO-stability  implies $H_\infty$-stability, and for finite-dimensional
LTI-systems the two are equivalent. There exist infinite dimensional LTI-systems which are $H_\infty$-stable,
but not BIBO stable. However, 
if the system is $H_\infty$-stable and
exponentially stabilizable and detectable, then it is exponential stable \cite{morris}, and that
implies BIBO stability. The latter because if the growth rate of $A$ is $<0$, then there exists $a > 0$ 
such that $|ce^{At}b| \leq Me^{-at}$, which implies integrability of $ce^{At}b$.
\end{remark}

\noindent
{\bf Results.} 
Experiments with different tolerances $\tau = 1.0$ and $\tau = 0.1$ were performed.
In each case the parameter $c$ varied in the interval $[-6.2,-0.5]$ and optimization led to an asymptotic sector, 
see Fig. \ref{fig4}. Smaller values of $c$ lead to larger aperture in the sectors.
Two scenarios were selected and underwent non-linear simulations
with three types of disturbances shown in Figs. \ref{fig7} and \ref{fig8}. The resulting asymptotic sectors are shown in Fig. \ref{as_sectors},
and typical optimized impulse responses are shown in Fig. \ref{responses}.

\section{Extension to multi-dimensional non-linearity}
\label{mimo}

In order to extend our algorithm to systems (\ref{NL1}) with
multi-dimensional non-linearity,
we
consider a feedback loop between an LTI-system $G$ and the non-linearity $\Delta$:
\begin{align}
\label{in_out}
    \begin{split}
        v &= Gw +  f\\
        w&= \Delta(v)+e
    \end{split}
\end{align}
as shown in
Fig. \ref{fig_pos_neg} (left). 
Well posedness of  (\ref{in_out})  in the $L_2$-sense 
means that $G,\Delta$ are $L_2$-bounded causal operators on $L_{2e}$, and that 
the map $(v,w)\to (e,f)$ has a causal inverse on $L_{2e}$. The system is $L_2$-stable
if this inverse is bounded, i.e., if there exists a constant $c>0$ with
$$
|v|_2^2 + |w|_2^2 \leq c\left(|f|_2^2+|e|_2^2 \right)
$$
for any solution of (\ref{in_out}). 
Since we are interested in BIBO-stability, we also need the corresponding notions in the time-domain
$L_\infty$-sense.

\begin{definition}{\rm
The feedback connection $(G,\Delta)$ is well-posed in the time-domain $L_\infty$ sense if
$\Delta$ and $G$ both map $L_{\infty e}$ into $L_{\infty e}$, and if the map $(v,w) \to (e,f)$ from (\ref{in_out})
has a causal inverse  $(e,f) \to (v,w)$ on the spaces
$L_{\infty e} \times L_{\infty e} \to L_{\infty e} \times L_{\infty e}$. } 
\end{definition}

\begin{definition}
{\rm
The $L_\infty$-well-posed feedback connection $(G,\Delta)$ is BIBO-stable, if in the setting of (\ref{in_out}) there exist
constants $k_1>0$, $k_2>0$, such that
\[
|v|_\infty + |w|_\infty \leq k_1\left( |e|_\infty + |f|_\infty \right)+k_2
\]
for all $e,f\in L_\infty([0,\infty),\mathbb R^n)$.} 
\end{definition}

We now investigate ways in which the steps of the algorithm in Section \ref{algorithm} may be extended to MIMO
non-linearity (\ref{in_out}).

\subsection{Extending the mixed $H_\infty/H_\infty$-program}
\label{4.1}
Extending step 4 to MIMO non-linearity 
leads to Integral Quadratic Constraints (IQC), where
$\Delta$ and $G$ in loop as in (\ref{in_out}) satisfy the 
quadratic constraints induced by a multiplier $\Pi=\Pi^\sim$:
\begin{equation}
\label{both_iqcs}
\left\langle \begin{bmatrix} u\\\Delta(u)\end{bmatrix} , \Pi \begin{bmatrix} u\\\Delta(u)\end{bmatrix}\right\rangle_T \leq 0, \quad \left\langle \begin{bmatrix} Gu\\u\end{bmatrix} , \Pi \begin{bmatrix} Gu\\u\end{bmatrix}\right\rangle_T \geq 0
\end{equation}
for every $u\in L_{2e}$ and every $T\geq 0$. While \cite{rantzer} 
assures $L_2$-stability of the loop if one of the inequalities is satisfied strictly, the
crucial question is how the IQC for $G$ may be verified
algorithmically.
In the literature these are traditionally transformed to LMIs, but
in synthesis lead to BMIs, which are known to encounter numerical difficulties. This was
recognized in \cite{iqc,prot,mu}, where non-differentiable optimization techniques
in tandem with Hamiltonian tests for function evaluations \cite{BB:92} were preferred instead. Recently this line
has
been further perfected in \cite{alberto1,pascal}.

In \cite{pascal}, the authors
obtain a mixed $H_\infty$/$H_\infty$-program expanding on (\ref{second}) for
J-spectral factorable multipliers 
$\Pi(s) = \Psi^\sim(s) J \Psi(s)$, where $J = [I_p , 0;0,-I_m]$
and $\Psi(s)$ is a $L_2$-bistable rational system. 
Defining processes
\begin{equation}
    \label{well}
\begin{bmatrix} \Psi_{11} u + \Psi_{12} \Delta(u)\\
\Psi_{21} u + \Psi_{22}\Delta(u)\end{bmatrix} =: \begin{bmatrix} \widetilde{u}\\\widetilde{\Delta}(\widetilde{u})\end{bmatrix},
\quad
\begin{bmatrix} \Psi_{11} Gu + \Psi_{12}u \\\Psi_{21}Gu + \Psi_{22}u\end{bmatrix}
=: \begin{bmatrix} \widetilde{G}\widetilde{u} \\ \widetilde{u}\end{bmatrix},
\end{equation}
$L_2$-stability of the loop $(G,\Delta)$ is equivalent to $L_2$-stability
of the loop $(\widetilde{G},\widetilde{\Delta})$. Assuming  that $\Delta$
is square, it follows from \cite[Thm. 5]{pascal} that 
$\widetilde{G} = (\Psi_{11}G+\Psi_{12})(\Psi_{21}G+\Psi_{22})^{-1}$
is well-posed, and the IQC for $G$ in (\ref{both_iqcs}) is transformed into
$\|\widetilde{G}\|_\infty \leq 1$. Similarly, with $\circ$ denoting map or relation composition,
$\widetilde{\Delta} = (\Psi_{22}\Delta + \Psi_{21}) \circ (\Psi_{12}\Delta + \Psi_{11})^{-1}$,
and due to (\ref{both_iqcs}), left,  this process  is
an $L_2$ contraction, i.e., signals $z_1= \Psi_{11} u + \Psi_{12}\Delta(u)$ and $z_2=\Psi_{21}u+\Psi_{22}\Delta(u)$
satisfy $\int_0^T |z_2(t)|^2 dt \leq \int_0^T |z_1(t)|^2 dt$ for all $T > 0$. 
What is not clear is whether $\widetilde{\Delta}$ is a mapping, because the argument which worked for
$\widetilde{G}$ in \cite{pascal} hinges on linearity.

We now present an alternative way to obtain a mixed $H_\infty/H_\infty$-program,
which gives an explicit loop transformation and, as we shall see, is also applicable to positivity type
factorizations.
Consider again IQC multipliers $\Pi(s)$ factored as
\begin{equation}
\label{eq_factor}
\Pi(j\omega) = \Psi(j\omega)^* P \Psi(j\omega)
\end{equation}
for a bistable LTI-system $\Psi(s)$ and a static invertible $P=P^T$. 
Such factorizations exist for rational $\Pi=\Pi^\sim$  
if $\Pi$ has neither poles nor zeros on $j\mathbb R$
and allows  no equalizing vectors, i.e., no $u\in \mathscr H_2$, $u\not=0$, with $\Pi u \in \mathscr H_2^\perp$;
cf. \cite{meinsma}. 
In particular, positive-negative multipliers satisfying
$\Pi_{11}(j\omega) \succeq \epsilon I$ and $\Pi_{22}(j \omega) \preceq -\epsilon I$ for some $\epsilon >0$
admit such factorizations \cite{seiler_2015}. 
There is no loss of generality
in assuming that both IQCs in (\ref{both_iqcs}) are satisfied strictly. 
Now define a new augmented interconnection $({G}_a,{\Delta}_a)$ as
\begin{equation}
    \label{Delta_a}
{G}_a = \Psi \begin{bmatrix} -I&2G\\0&I\end{bmatrix} \Psi^{-1}, \quad
{\Delta}_a = \Psi \circ \begin{bmatrix} I&0\\2\Delta&-I\end{bmatrix} \circ \Psi^{-1},
\end{equation}
then by \cite[Thm. 2 (1)]{scherer}, $L_2$-stability of $(G,\Delta)$ is equivalent
to $L_2$-stability of $({G}_a,{\Delta}_a)$. 

Adopting $\Pi_{11} \succeq \epsilon I$
and $\Pi_{22} \preceq -\epsilon I$ for some $\epsilon >0$, it follows from \cite[Thm. 2 (2)]{scherer} that 
${\Delta}_a$, $G_a$ satisfy IQCs for the passivity multiplier
${P}_a=[0, P;P,0]$ strictly, i.e., 
\begin{align}
\label{He}
\begin{split}
\int_0^T {p}_a(t)^T {G}_a(P^{-1} {p}_a)(t)  dt &\leq -\epsilon \int_0^T {p}_a(t)^T {p}_a(t)  dt\\  
\int_0^T {p}_a(t)^T  P {\Delta}_a({p}_a)(t) dt &\geq \epsilon \int_0^T {p}_a(t)^T p_a(t) dt 
\end{split}
\end{align}
for some $\epsilon > 0$ and every $T\geq 0$, where $p_a=(p,q)^T$.

The
inequalities in (\ref{He}) are now turned into bounded gain conditions using M\"obius or bilinear transformations. 
We introduce 
\begin{equation}
    \label{Ge}
{{G}_e}:= \mathcal B \star {G}_a P^{-1} = (G_aP^{-1}-I)^{-1}(G_aP^{-1}+I)
\end{equation}
and
\begin{equation}
    \label{Delta_e}
{\Delta}_e = P\circ {\Delta}_a  \star \mathcal B = (I+P\circ \Delta_a)^{-1} \circ (I-P\circ \Delta_a)
\end{equation}
with $$\mathcal B:= \begin{bmatrix} -I & \sqrt{2} I\\ -\sqrt{2}I & I\end{bmatrix},$$
where $\star$ is the Redheffer star product \cite{redheffer}. 
Here $(I+P\Delta_a)^{-1}$ and $(G_aP^{-1}-I)^{-1}$ are well-defined and stable due to
(\ref{He}) and
the passivity theorem, hence $\Delta_e,G_e$ are well-defined and $L_2$-stable.
Indeed, $(I+P\circ \Delta_a)^{-1}$ is the negative feedback loop between
the upper block $I$ and the lower block $P\circ \Delta_a$. Since $P\circ \Delta_a$ is strictly passive by (\ref{He}) and $I$ is passive, stability follows from the passivity theorem
\cite{zames:66}. A similar argument
applies to $(G_aP^{-1}-I)^{-1}$.
Owing to
$\mathcal B\star \mathcal B = I^\sharp = [0,I;I,0]$, the unit of the star
product, we get (stability) loop invariance
$(G_a,\Delta_a)\cong (G_aP^{-1},P\circ \Delta_a) \cong (\mathcal B\star G_aP^{-1}, P\circ \Delta_a\star \mathcal B)=
(G_e,\Delta_e)$.
This means the passivity-type conditions (\ref{He}) 
are equivalent to
bounded-gain conditions 
\begin{align}
\label{Bounded}
\begin{split}
\int_0^T  \|{G}_e({p}_e)(t)\|^2  dt &\leq (1-\epsilon) \int_0^T \|{p}_e(t)\|^2  dt  \\
\int_0^T \|  {\Delta}_e({p}_e)(t)\|^2 dt &\leq (1-\epsilon) \int_0^T \|{p}_e(t)\|^2  dt \,,
\end{split}
\end{align}
for some $\epsilon >0$ and every $T\geq 0$, where $p_e=(p,q)^T$. See
\cite[pp. 215-16]{desoer} for a proof, which also applies to  the non-linear case.

From $\Delta_e=P\circ \Delta_a\star \mathcal B$ we have $\Delta_e\star \mathcal B = P\circ\Delta_a\star\mathcal B\star \mathcal B = P\circ \Delta_a\star I^\sharp = P\circ \Delta_a$,
hence
$\Delta_a = P^{-1} \circ \Delta_e\star \mathcal B$. That gives $\begin{bmatrix} I&0\\2\Delta &-I\end{bmatrix}
= \Psi^{-1} P^{-1}\circ \Delta_e\star \mathcal B \Psi$. Hence,
\begin{equation}
    \label{all_Deltas}
\Delta = \begin{bmatrix}  \frac{I}{\sqrt{2}}\\\frac{I}{\sqrt{2}}\end{bmatrix}^T \Psi^{-1} (P^{-1} \circ \Delta_e\star \mathcal B) \Psi \begin{bmatrix} \frac{I}{\sqrt{2}}\\\frac{I}{\sqrt{2}}\end{bmatrix},
\end{equation}
which gives the inverse operation to $\Delta \to \Delta_e$
in (\ref{Delta_e}).
What we have obtained is a parametrization of all non-linearities $\Delta$ derived from
$L_2$-contractions $\Delta_e$ via the loop transformation through $\Psi(s), P$, or equivalently,
all non-linearities satisfying IQCs with factorable multiplies $\Pi = \Psi^*P\Psi$. For these $\Delta$
the IQC-stability theorem can now be reduced to the small gain theorem \cite{zames:66}.

\begin{theorem}
\label{theorem2}
{\rm\bf  (IQC as $H_\infty$ constraint).}
Suppose $(G,\Delta)$ is loop transformed to $(G_e,\Delta_e)$, where $\Delta$ satisfies the
IQC with multiplier
$\Pi=\Psi^*P\Psi$ factored with bistable  $\Psi(s)$ and invertible
$P=P^T$. Then $\|G_e\|_\infty < 1$ implies $L_2$-stability of the loop $(G,\Delta)$.
\hfill $\square$
\end{theorem}

Applying the loop transformation $(G,\Delta) \cong (G_e,\Delta_e)$ to the closed loop system
$\mathcal F_l(G,K)$ leads to $(\mathcal F_l(G,K),\Delta) \cong (\mathcal F_l(G,K)_e,\Delta_e)$.
This allows us now to extend step 4 of the algorithm to IQCs.

\begin{corollary}
\label{theorem3}
Suppose the mixed $H_\infty/H_\infty$-synthesis program
\begin{eqnarray}
\label{new_first}
\begin{array}{ll}
\mbox{\rm minimize} & \| \mathcal F_l(G,K)_e\|_\infty \\
\mbox{\rm subject to} & \|T_{wz}(G,K)\|_\infty \leq (1+\tau)\gamma_\infty\\
&K \in \mathscr K
\end{array}
\end{eqnarray}
admits an optimal solution $K^\sharp\in \mathscr K$ satisfying $\|\mathcal F_l(G,K^\sharp)_e\|_\infty < 1$.
Then $K^\sharp$ stabilizes the loop $(G,\Delta)$ in the $L_2$-sense, and linearized closed loop performance
is degraded over
nominal performance $\gamma_\infty$ by no more than
the factor $1+\tau$.
\hfill $\square$
\end{corollary}

\begin{remark}
Program (\ref{new_first}) is now a natural MIMO-extension of (\ref{first}).
It can be efficiently solved by the method of
\cite{AN06a,an_disk:05} available in the {\tt systune} package of \cite{CST2020b,apkarian2014multi,an_disk:05}.
This is numerically 
preferable to transforming IQCs to  BMIs. 
With the recent extension of non-smooth $H_\infty$-synthesis
in \cite{AN:18,apkarian:19,an1}
it becomes even possible to address (\ref{new_first}) for
infinite-dimensional systems with infinite-dimensional multipliers $\Psi(s)$.
\end{remark}

\begin{remark}
An advantage of this construction is that when $(G,\Delta)$ satisfies
an IQC with positivity multiplier as in (\ref{He}), then going from $(G,\Delta)$
to $(G_a,\Delta_a)$ can be skipped and we build
$(G_e,\Delta_e)$ directly without the augmentation (\ref{Delta_a}). 
\end{remark}

For multipliers
$\Pi(s)$ with lower triangular factorizations both approaches
(\ref{well}) and the augmentation $G\to G_a\to G_e$ lead to the same result. Suppose
\begin{equation}
    \label{triag_factor}
\Pi(j\omega)= \Psi^T(-j\omega) P\Psi(j\omega), \quad
 \Psi = \begin{bmatrix} \Psi_{11}&0\\\Psi_{21} & \Psi_{22} \end{bmatrix} \quad
P=\begin{bmatrix}
I&0\\0&-I
\end{bmatrix},
\end{equation}
with $\Psi_{11},\Psi_{11}^{-1},\Psi_{21},\Psi_{22},\Psi_{22}^{-1}$ stable.
Then the transformed non-linear operator and LTI-system
\begin{equation}
    \label{tilde_transform}
\widetilde{\Delta} = \Psi_{21} \Psi_{11}^{-1} + \Psi_{22} \circ \Delta \circ \Psi_{11}^{-1},\qquad
\widetilde{G}=\Psi_{11}G(\Psi_{22}+\Psi_{21}G)^{-1},
\end{equation}
give an equivalent loop $(G,\Delta)\cong (\widetilde{G},\widetilde{\Delta})$, where
the IQC is transformed to a Small-Gain condition $|\widetilde{\Delta}(\widetilde{v})|_2 \leq |\widetilde{v}|_2$,
$\|\widetilde{G}\|_\infty < 1$, now with $\widetilde{\Delta}$ and $\widetilde{G}$ of the same dimension
as $\Delta,G$. Since  $\|G_e\|_\infty < 1$ is equivalent to $\|\widetilde{G}\|_\infty < 1$ in the case
(\ref{triag_factor}), 
Theorem \ref{theorem2} implies:

\begin{corollary}
\label{theorem4}
{\rm \bf (Triangular transform)}.
Suppose a non-linearity $\Delta$ can be loop transformed to a $L_2$-contraction
$\widetilde{\Delta}$ in {\rm (\ref{tilde_transform})} using a lower triangular factorization 
{\rm (\ref{triag_factor})}. %
Suppose the transformed {\rm LTI}-system $\widetilde{G}=\Psi_{11} G \left( I+\Psi_{22}^{-1} \Psi_{21} G \right)^{-1} \Psi_{22}^{-1}$ is stable and satisfies $\|\widetilde{G}\|_\infty < 1$. Then the loop $(G,\Delta)$
is stable in the $L_2$-sense.
\end{corollary}

\begin{proof}
Indeed, with (\ref{triag_factor}) the expression for $G_e$ simplifies to
$$
G_e = \begin{bmatrix} 0 & \Psi_{11} G \left(\Psi_{22}+\Psi_{21} G  \right)^{-1}\\ \Psi_{21}\Psi_{11}^{-1} & 0 \end{bmatrix}
$$
with $\widetilde{G}$ in the upper right corner. 
Since $\Pi_{11} = \Psi_{11}^*\Psi_{11} - \Psi_{21}^*\Psi_{21} \succ 0$ by hypothesis, we have
$\|\Psi_{21}\Psi_{11}^{-1} \|_\infty < 1$, hence $\|G_e\|_\infty < 1$ is
equivalent to $\|\widetilde{G}\|_\infty < 1$ as claimed.
\hfill $\square$
\end{proof}

This 
may also be seen from (\ref{well}); see also Fig. \ref{fig_pos_neg}.  For
upper triangular $\Psi$ we obtain the analogous result with
\begin{equation}
\label{upper_triangular}
\widetilde{\Delta} = \Psi_{22} \Delta \circ (\Psi_{12}\Delta + \Psi_{11})^{-1},
\quad \widetilde{G}=(\Psi_{11}G+\Psi_{12}) \Psi_{22}^{-1}.
\end{equation}
These results no longer require any reference to IQCs or multipliers.

\subsection{Extending the mixed peak-gain/$H_\infty$-program}
\label{4.2}
The one-dimensional peak-gain norm \cite{boyd_barratt} allows several extensions  to MIMO systems,
because we can replace the absolute value $|x|$, $x\in \mathbb R$,  by any of the equivalent vector norms in $\mathbb R^n$. If we define 
a signal norm
on $L_\infty([0,\infty),\mathbb R^n)$ by
$$
|x|_{\infty,p} = \sup_{t\geq 0} |x(t)|_p
$$
with $|v|_p$ the $p$-norm of $v\in \mathbb R^n$, $1 \leq p \leq \infty$, then with the notation adopted from
\cite{induced_norms} {\it any} induced system norm
\begin{equation}
    \label{norms}
\|G\|_{(\infty,p),(\infty,q)} = \sup_{x \not=0} \frac{|G\ast x|_{\infty,q}}{|x|_{\infty,p}}
\end{equation}
is a valid MIMO extension of $\|\cdot\|_{\rm pk\_gn}$. 
The peak-gain  norm, to which we give preference here, is the special case
$\|G\|_{\rm pk\_gn} = \|G\|_{(\infty,\infty),(\infty,\infty)}$, 
but all norms (\ref{norms}) are equivalent. 

\begin{theorem}
\label{small_gain}
Suppose the non-linear operator in {\rm (\ref{in_out})} satisfies $|\Delta(t,x)|_{p} \leq |x|_{q}$ 
for every $|x|_{q} > M$ and $|\Delta(t,x)|_{p}\leq L$ for every $|x|_{q}\leq M$. If the LTI-system
$G$ satisfies $\|G\|_{(\infty,p),(\infty,q)} < 1$, then the closed loop is
BIBO-stable with
$|v|_{\infty,q} + |w|_{\infty,p} \leq k_1 (|e|_{\infty,p}+|f|_{\infty,q}) + k_2$ for all
$e,f \in L_{\infty e}$.
\end{theorem}

\begin{proof}
Put $\Delta_1(t,x) = \Delta(t,x) \chi_{\{|x|_{q} \leq M\}}(x)$, $\Delta_2(t,x)=\Delta(t,x) \chi_{\{|x|_{q} > M\}}(x)$,
then we have $\sup_{t\geq 0} |\Delta_1(t,v(t))|_p \leq L$, while $\sup_{t\geq 0} |\Delta_2(t,v(t))|_p \leq \sup_{t\geq 0}|v(t)|_q$ by hypothesis. Hence, assuming $|Gx|_{\infty,q} \leq (1-\delta) |x|_{\infty,p}$ for some $0 < \delta < 1$,
\begin{align*}
    |v|_{\infty,q} &\leq |w'|_{\infty,q} + |f|_{\infty,q} \\
    &\leq |G \Delta_1(v)|_{\infty,q} + |G\Delta_2(v)|_{\infty,q} + |Ge|_{\infty, q} + |f|_{\infty,q} \\
    &\leq \|G\|_{(\infty,p)(\infty,q)} \left( |\Delta_1(v)|_{\infty,p} + |\Delta_2(v)|_{\infty,p} + 
    |e|_{\infty,p}\right) + |f|_{\infty,q}\\
    &\leq (1-\delta)L + (1-\delta) |v|_{\infty,q} + (1-\delta)|e|_{\infty,p} + |f|_{\infty,q}
\end{align*}
hence
$$
\delta |v|_{\infty,q} \leq (1-\delta) |e|_{\infty,p} + |f|_{\infty,q} + (1-\delta) L.
$$
On the other hand
$$
|w|_{\infty,p} \leq |\Delta_1(v)|_{\infty,p} + |\Delta_2(v)|_{\infty,p} + |e|_{\infty,p}
\leq L + |v|_{\infty,q} + |e|_{\infty,p}.
$$
Combining the two implies the estimate.
\hfill $\square$
\end{proof}

For this result compare the more general \cite{mareels_hill:92}, \cite{teel}.
The above proof is standard and included for convenience.
This gives us now a clue how to extend asymptotic constraints as encountered in
Section \ref{asymptotic}
to MIMO non-linearities.

\begin{definition}
\label{static}
{\bf (Asymptotic $L_\infty$-contraction)}. {\rm 
A non-linear operator $\Delta:[0,\infty) \times \mathbb R^n\to  \mathbb R^n$ is called an asymptotic 
$L_\infty$-contraction
if there
exist $L,M > 0$ such that
$|\Delta(t,x)|_\infty \leq |x|_\infty$ for all $|x|_\infty > M$, $t\geq 0$, and
$|\Delta(t,x)|_\infty \leq L$ for all $|x|_\infty \leq M, t\geq 0$.
}
\end{definition}

\begin{remark}
The proof of Theorem \ref{small_gain} shows that an asymptotic  $L_\infty$-contraction
satisfies $|\Delta(t,x)|_\infty \leq |x|_\infty + k$ for all $x$. Conversely, suppose
we have $|\Delta(t,x)|_\infty \leq |x|_\infty + k$ for all $x$.
Then for every $\epsilon >0$ there exists $M>0$ such that $|\Delta(t,x)|_\infty < (1+\epsilon)|x|_\infty$ for all $|x|_\infty > M$.
For suppose on the contrary that there exist $x_n$ with $|x_n|_\infty \to \infty$ such that 
for some $\epsilon >0$
$|\Delta(t,x_n)|_\infty \geq (1+\epsilon)|x_n|_\infty$, 
then $1+\epsilon \leq |\Delta(t,x_n)|_\infty/|x_n|_\infty \leq 1 + k/|x_n|_\infty\to 1$,
a contradiction. Since for the LTI-system we request the strict inequality $\|G\|_{\rm pk\_gn} < 1$,
both conditions for $\Delta$ may be used indifferently in the small gain theorem.
\end{remark}

\begin{figure}[ht!]
\includegraphics[scale=0.9]{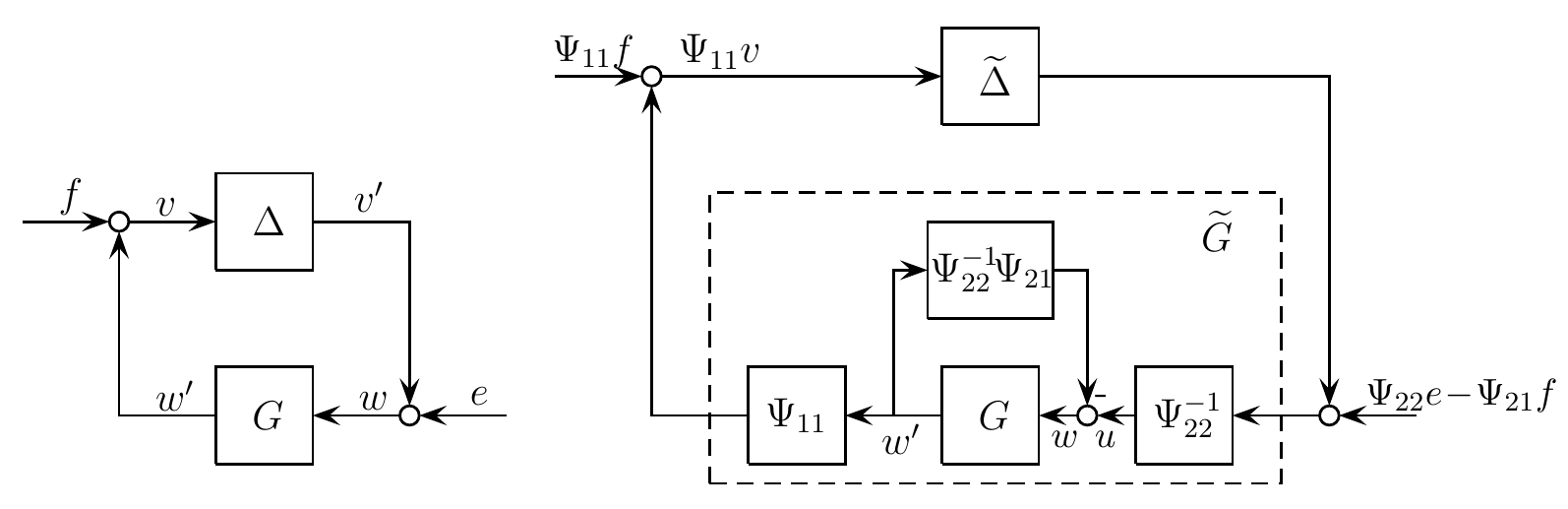}
\caption{Loop-transformation $(G,\Delta)$ to $(\widetilde{G},\widetilde{\Delta})$. 
\label{fig_pos_neg}}
\end{figure}

\begin{corollary}
{\bf (Triangular transform)}.
Let $(G,\Delta)$  be $L_\infty$ well-posed, and suppose
$\Delta$ can be loop-transformed via a lower triangular $L_\infty$-bistable $\Psi$ to an asymptotic $L_\infty$-contraction $\widetilde{\Delta}$. 
Suppose $\widetilde{G}=\Psi_{11}G(\Psi_{22}+\Psi_{21}G)^{-1}$
is
$L_\infty$ well-posed and satisfies $\|\widetilde{G}\|_{\rm pk\_gn} < 1$.
Then the loop {\rm (\ref{in_out})} is BIBO-stable.
\end{corollary}

\begin{proof}
This refers to (\ref{tilde_transform})  shown in
Fig. \ref{fig_pos_neg}, where $(G,\Delta)$ is loop transformed to
$(\widetilde{G},\widetilde{\Delta})$ in such a way that BIBO-stability of
$(G,\Delta)$ is equivalent to BIBO-stability of $(\widetilde{G},\widetilde{\Delta})$.
But $(\widetilde{G},\widetilde{\Delta})$ is amenable to Theorem \ref{small_gain}, hence
$\|\widetilde{G}\|_{\rm pk\_gn} < 1$ implies BIBO-stability of the loop.
\hfill $\square$
\end{proof}

We now extend program (\ref{second}) to the MIMO-case.
In the case of Fig. \ref{optimistic},  we transform the non-linearity
$\Delta$ to an asymptotic contraction $\widetilde{\Delta}$ via (\ref{tilde_transform}). Now consider a plant
$$P:\begin{bmatrix} q\\z\\y\end{bmatrix} = \begin{bmatrix} P_{11} & 0 & P_{13}\\0&P_{22}&P_{23}\\P_{31}&P_{32}&P_{33}\end{bmatrix} \begin{bmatrix} p\\w\\u \end{bmatrix},
P_1=\begin{bmatrix} P_{11}&P_{13}\\P_{31}&P_{33}\end{bmatrix},
P_2 = \begin{bmatrix} P_{22}&P_{23}\\P_{32}&P_{33}\end{bmatrix}$$
then the extension of (\ref{second}) has the form of the mixed program 
\begin{align*}
&(\ref{second}'') \hspace{4cm}
\begin{array}{ll}
\mbox{minimize} & \| \mathcal F_l(P_1,K)^\sim \|_{\rm pk\_gn}\\
\mbox{subject to}& \| \mathcal F_l(P_2,K)\|_\infty \leq (1+\tau)\gamma_\infty\\
&K \in \mathscr K
\end{array}
\hspace*{4cm}
\end{align*}
where $\sim$ indicates that the loop transformation $\widetilde{G}$ of Fig. \ref{fig_pos_neg} 
is applied to the controlled system
$\mathcal F_l(P_1,K)$. The program is successful as soon as a structured LTI 
controller $K^*\in \mathscr K$ is found which stabilizes $P_1$ in the BIBO-sense, stabilizes
$P_2$ exponentially, and achieves $\|\mathcal F_l(P_1,K^*)^\sim\|_{\rm pk\_gn} < 1$.

\begin{remark}
We could also use the general loop transformation
(\ref{all_Deltas}). Suppose $\Delta$ is obtained from a $|\cdot|_\infty$ contraction
$\Delta_e$ via (\ref{all_Deltas}), where $\Delta_e \star \mathcal B$ is $L_\infty$-well-posed. Then
a sufficient condition for BIBO-stability of the loop is $L_\infty$-well-posedness of $G_e$ in 
(\ref{Ge}) with $\|G_e\|_{\rm pk\_gn}<1$, and this includes both cases (\ref{tilde_transform}), (\ref{upper_triangular}).
\end{remark}

\subsection{Asymptotic $L_\infty$-contractions}
\label{sect_asymptotic_mimo}
In this section, we collect a variety of examples of MIMO
non-linearities which may be assessed by way of asymptotic $L_\infty$-contractions.

\begin{example}
\label{example1}
Consider a non-linearity $\Delta(t,q)$ in loop $(G,\Delta)$ as
\begin{align*}
G: \quad \begin{split}
\dot{x} &= Ax + Bp\\
q&= Cx
\end{split}
\qquad \qquad p(t)=\Delta(t,q(t)).
\end{align*}
We say that $\Delta$ is 
asymptotically polyhedral, if there
exist polyhedral norms $|\cdot|_\square$ and $|\cdot|_\triangle$ on $\mathbb R^n$ and $k>0$ such that
$|\Delta(t,q)|_\square \leq |q|_\triangle + k$ for all $q\in \mathbb R^n$ and all $t\geq 0$.  Now polyhedral norms
are of the form
$$
|x|_\square = \sup_{i=1,\dots,m} \left| \sum_{j=1}^n \tau_{ij} x_j\right| = |Tx|_\infty, \qquad
|y|_\triangle = \sup_{k=1,\dots,p} \left| \sum_{j=1}^n \sigma_{kj} y_j\right| = |Sy|_\infty
$$
for certain $T\in \mathbb R^{m\times n}$, $S\in \mathbb R^{p\times n}$ with $\{x\in \mathbb R^n: |Tx|_\infty \leq 1\}, \{y\in \mathbb R^n:  |Sy|_\infty \leq 1\}$ bounded.
The latter means $T,S$ are injective.  Let $T^+,S^+$ be left inverses, $T^+Tp=p$, $S^+Sq=q$. 
Then we have
$$
|T \Delta(t,S^+Sq)|_\infty =
| T \Delta(t,q) |_\infty \leq |Sq|_\infty + k
$$
so on introducing the new variables $q_1= Sq$, $p_1=Tp$, we have
a new non-linearity $\widetilde{\Delta}(t,\cdot) = T \circ \Delta(t,\cdot) \circ S^+$ which satisfies
$|\widetilde{\Delta}(t,q_1)|_\infty \leq |q_1|_\infty + k$ for all $q_1$ and $t\geq 0$.

The non-linearity $\Delta$ being in loop with $G$, 
we transform this to bring $\widetilde{\Delta}$ in loop with:
\begin{align*}
\widetilde{G}: \quad \begin{split}
\dot{x} &= Ax + BT^+p_1\\
q_1&= SCx\\
\end{split}
\end{align*}
This means the non-linear loop $(G,\Delta)$ is BIBO-stable if $\|\widetilde{G}\|_{\rm pk\_gn} < 1$.
This is a special case of the transform (\ref{tilde_transform}).
\end{example}

\begin{example}
\label{example2}
The following is a concretization.
We call $\Delta$ differentiable at infinity if there exists a matrix 
$\Delta_\infty: = \Delta'(\infty)\in \mathbb R^{n\times n}$ such that
$$
\lim_{|x|\to \infty} \frac{|\Delta(t,x) - \Delta'(\infty)x|}{|x|} = 0
$$
uniformly over $t\geq 0$.
Here we may choose arbitrary norms in numerator and denominator. Consider the case where $\Delta'(\infty)\not=0$
and choose regular matrices
$T,S\in \mathbb R^{n\times n}$ such that $T\Delta'(\infty)S^{-1} = {\rm diag}(1-\epsilon,\dots,1-\epsilon,0,\dots,0)=:J_\epsilon$, where the diagonal
has rank$(\Delta'(\infty))$ many entries $1-\epsilon$. Now choose the norms $|y|_\square =|Ty|_\infty$ and $|x|_\triangle = |Sx|_\infty$,
then
\begin{align*}
\frac{|\Delta(t,x)|_\square}{|x|_\triangle}& \leq \frac{| \Delta(t,x) - \Delta'(\infty)x|_\square}{|x|_\triangle} +
\frac{|\Delta'(\infty)x|_\square}{|x|_\triangle}\\
&=o(1) + \frac{|T \Delta'(\infty) S^{-1}Sx|_\infty}{|Sx|_\infty}\\
&\leq o(1)+ \frac{\normi{J_\epsilon}_\infty |Sx|_\infty}{|Sx|_\infty} \leq o(1) + 1-\epsilon,
\end{align*}
where $\normi{J_\epsilon}_\infty=1-\epsilon$ is the maximum row sum norm.
Choosing $M>0$ such that $|o(1)| < \epsilon/2$ for $|x|_\triangle > M$, we arrive at
$|\Delta(t,x)|_\square \leq (1-\epsilon/2) |x|_\triangle + \sup\{|\Delta(t,x)|_\square: |x|_\triangle \leq M\}=:
(1-\epsilon/2) |x|_\triangle +k$. This means every non-linearity which is differentiable at infinity admits
asymptotic $L_\infty$-constraints.
\end{example}

\begin{remark}
Suppose a non-linearity $\Delta:[0,\infty)\times \mathbb R^n \to \mathbb R^n$ satisfies 
$|\Delta(t,q)|_2 \leq |q|_2 + k$ for some $k\geq 0$ and all $q\in \mathbb R^n$, $t\geq 0$. 
Then in Theorem \ref{small_gain} we would prefer
the system norm $\|G\|_{(\infty,2),(\infty,2)}$. Unfortunately, no computable expression is 
currently known for
this norm, so its optimization is presently impossible.
\end{remark}

\begin{remark}
\label{makeshift}
For the case $|\Delta(t,q)|_2 \leq |q|_2+k$ we have the following makeshift alternative.
Choose approximations
$P_1 \subset B(0,1) \subset P_2$ by polytopes $P_1,P_2$, then $\Delta(P_1) \subset P_2$
asymptotically, so we are in the situation of Example \ref{example1} and we may work with
$\|\cdot\|_{\rm pk\_gn}$. In the case
$|\Delta(t,q)|_2 \leq (1-\epsilon)|q|_2 + k$ we may even obtain this with $P_1=P_2
\subset \{q: |q|_2\leq 1\}$ and vertices on $\{q: |q|_2=1\}$.
\end{remark}

\begin{example}
\label{example3}
MIMO sectors are defined via 
symmetric matrices $A,B \in \mathbb S^n$ satisfying $A \prec B$. 
A mapping $\phi:\mathbb R^n \to \mathbb R^n$ with $\phi(0)=0$
is in the sector {\bf sect}$(A,B)$, noted $\phi \in {\bf sect}(A,B)$, if $(\phi(x) - Ax)^T(\phi(x)-Bx) \leq 0$
for all $x\in \mathbb R^n$.

With the choice $R = \frac{1}{2} (B-A) \succ 0$ and $C = \frac{1}{2} (B+A)$ we find
that $\psi(x) = \phi(x) - Cx$ satisfies $\psi(x)^T\psi(x) \leq x^TR^TRx$, so we get a norm
bound $|\psi(x)|_2 \leq |Rx|_2$, and if we define
$\Delta = \psi\circ R^{-1}=(\phi-C)\circ R^{-1}$, then $|\Delta(y)|_2 \leq |y|_2$. 
When we allow $|\Delta(y)|_2 \leq |y|_2+k$ for some $k \geq 0$ and all $y$, this is
a typical application of 
the two previous remarks, where we would like to apply Theorem \ref{small_gain} with $\|\cdot\|_{(\infty,2),(\infty,2)}$.
\end{example}

\begin{example}
\label{example4}
As a concretization
\cite{syazreen,bemporad}
consider a non-linearity $\phi$ generated by a convex quadratic program:
\begin{equation}
\label{QP}
\phi(x) = {\rm argmin} \{ \textstyle\frac{1}{2} v^THv - v^Tx: Lv \leq b\} 
\end{equation}
where $H \succ 0$ and $L\in \mathbb R^{m\times n}$, $b\geq 0, b\in \mathbb  R^m$ are fixed, and optimization is over $v\in \mathbb R^n$. 
Using the Kuhn-Tucker conditions
one verifies that $\phi$ satisfies the MIMO sector bound
$$
\phi(x)^T( H\phi(x)-x) \leq 0 \mbox{  for all $x$}.
$$
Then from the above
$\psi(x) = 2H^{1/2}\phi(H^{1/2}x) -x$
is a $|\cdot|_2$-contraction.
This means an asymptotic quadratic constraint for $\phi$ would lead
to $\|\cdot\|_{(\infty,2),(\infty,2)}$, which is, however, not available for computations.
A polyhedral approximation based on Remark \ref{makeshift} may be used instead.
\end{example}

\begin{example} (Continued).
\label{example5}
In the above case we can bring in $\|G\|_{\rm pk\_gn}$ directly, because the solution mapping
of a convex quadratic program with perturbation of the linear term or the constraints is known to
be piecewise affine \cite{bemporad,bank}, so it maps polyhedra to polyhedra. This allows
a construction as in Example \ref{example1}.

Since $H \succ 0$, (\ref{QP}) is equivalent to projecting $H^{-1}x$ orthogonally on the polyhedron
$\{v: Lv \leq b\}$ with regard to the Euclidean norm $|x|_H^2 =x^THx$. Let
$I \subset \{1,\dots,m\}$, $J=\{1,\dots,m\}\setminus I$, so that
$F=\{v\in \mathbb R^n: L_Iv=b_I, L_Jv\leq b_J\}$ is a face of the polyhedron, then projection of $x$
on $F$ is obtained as 
$$\phi(x) = H^{-1} \left( x-L_I^T(L_IH^{-1} L_I^T)^{-1} \left[ L_IH^{-1} x-b_I  \right]\right)$$ 
using generalized least squares. This shows that $\phi$ is piecewise affine. 
\end{example}

\begin{example}
\label{example6}
(Piecewise affine).
A non-linearity $\phi$ on $\mathbb R^n$ is {\it piecewise affine} if
there exist finitely many non-overlapping polyhedra $P_1,\dots,P_N$ 
with $\bigcup_{i=1}^N P_i=\mathbb R^n$ such that $\phi|P_i$ is affine, i.e.
there exist an affine mapping
$A_i(x) = b_i + L_ix$ with $A_i|P_i=\phi|P_i$. Here $L_i$ is the linear part
of $A_i$. For each polyhedron $P_i$
choose a Motzkin decomposition
$P_i=Q_i+C_i$ with $Q_i$ a polytope and $C_i$ a polyhedral cone.
Let $B_\triangle=\{x: |x|_\triangle \leq 1\}$ be the unit ball of a polyhedral norm, 
compute the polytopes $B_i = B_\triangle \cap C_i$ and
$B_i'=L_i(B_i)$. (Note that $L_i(C_i)$ is a polyhedral cone, so if it is bounded, it reduces
to $\{0\}$, in which case $B_i'=\{0\}$, too. Therefore only unbounded $L_i(C_i)$ have to be
considered).
Finally
let
$B'$ be the convex hull of $\bigcup_{i=1}^N B_i'$, then $B'$ is a  polytope containing $0$. 
Let $k_1=\max\{|b'|_\infty: b'\in B'\}$. 

We claim that there exists
a constant $k > 0$ such that
for every $x\in \mathbb R^n$, $\phi(x) \in |x|_\triangle B' + kB_\infty$. Indeed,
let $x\in P_i$, $x = y + z$ with $y\in Q_i$, $z\in C_i$. 
Then $\phi(x)=b_i+L_ix = b_i+L_iy + L_iz= b_i+L_iy + |z|_\triangle L_i(z/|z|_\triangle)
= b_i + L_iy_i + |z|_\triangle b'$ for some $b'\in B_i'$, using the fact that
$z/|z|_\triangle \in B_\triangle \cap C_i$, hence $L_i(z/|z|_\triangle)\in L_i(B_\triangle \cap C_i) = B_i'$. 

Now $|z|_\triangle = |x-y|_\triangle \leq |x|_\triangle + |y|_\triangle$, hence
$|z|_\triangle b' = (|x|_\triangle + |y|_\triangle) \rho b'$ for some $\rho\in [0,1]$,
and since $0\in B_i'$, we have $\rho b'=b''\in B_i'$. Then
$|z|_\triangle b' = |y|_\triangle b'' + |x|_\triangle b''$.

Altogether $\phi(x) = b_i + L_iy + |y|_\triangle b'' + |x|_\triangle b''$, and here
the term $b_i+L_iy + |y|_\triangle b''$ is bounded independently of $|x|_\triangle$, because
$y \in Q_i$ are bounded. We put 
$k_2 = \max_{i=1,\dots,N} |b_i|_\infty$, $k_3 = \max_{i=1,\dots,N} \max_{y\in Q_i} |L_iy|_\infty$
and $k_4=\max_{i=1,\dots,N} \max_{y\in Q_i} |y|_\infty$, then
$k = k_2 + k_3 + k_1k_4$.

Now let $B_\square = {\rm co}(B' \cup (-B'))$, then $B_\square$ is a symmetric polytope, hence its Minkowski functional is a 
polyhedral norm, $|\cdot|_\square$,
and we have shown $|\phi(x)|_\square \leq |x|_\triangle + k$
for all $x$. The  matrix $T$ defining $|\cdot|_\square$ can be obtained
from the polyhedral representation
$B'=\{x: Tx\leq 1\}$, which for moderate dimensions of $x$ can be
pre-computed. The matrix $S$ is obtained from $|x|_\triangle=|Sx|_\infty$.
\end{example}

\begin{example}
\label{example7}
We consider a numerical example with a non-linearity based on (\ref{QP}).
The parametric quadratic program
\begin{eqnarray*}
(P)_q
\begin{array}{ll}
\mbox{minimize} & \textstyle\frac{1}{2} x_1^2 + \textstyle\frac{1}{2}x_2^2 - q_1 x_1 - q_2x_2 \\
\mbox{subject to} & x_1-x_2 \leq 3\\
&x_1+x_2 \geq 0\\
&x_1\geq 0\\
&x_2\geq -1
\end{array}
\end{eqnarray*}
defines  a non-linear operator $\phi:\mathbb R^2\to \mathbb R^2$ via $p={\rm argmin}(P)_q$. Here $p=\phi(q)$ is
the orthogonal projection of $q$ on the polyhedron $P=\{x\in \mathbb R^2: Lx\leq b\}$, where
$L^T=\begin{bmatrix} 1&-1&-1&0\\-1&-1&0&-1\end{bmatrix}$, $b^T=\begin{bmatrix} 3&0&0&1\end{bmatrix}$.
For every face $F$ of $P$ the set $F'=\phi^{-1}(F)$ is a polyhedron and each $\phi|F':F'\to P$ is affine.
Here $P$ has 8 faces, three vertices, four facets, and $P$ itself. The unbounded faces are
$F_1=\{(2,-1)\} + \mathbb R_+(1,1)= Q_1+C_1$, $F_2=\mathbb R_+(0,1)= C_2$, and $P=Q+0^+P$, with
$0^+P=\{x: Lx\leq 0\}$. Now
$F_1' = \phi^{-1}(F_1) = \{(2,-1)\}+\mathbb R_+(1,1)+\mathbb R_+(1,-1)=Q_1+C_1'$ and
$\phi_1 = \phi|F_1'=A_1|F_1'$ for the affine mapping $A_1(x) = L_1x+b_1=\begin{bmatrix} {1}/{2}&{1}/{2}\\{1}/{2}&{1}/{2}\end{bmatrix} \begin{bmatrix}x_1\\x_2\end{bmatrix}+\begin{bmatrix}3/2\\-3/2\end{bmatrix}$. 
Similarly,
$F_2'= \phi^{-1}(F_2)= \mathbb R_+(-1,0)+\mathbb R_+(0,1)=C_2'$, and
$\phi_2 = \phi|F_2' = A_2|F_2'$ for the affine mapping
$L_2x=A_2(x) = \begin{bmatrix} 0&x_2\end{bmatrix}^T$, which is already linear. Clearly, $\phi^{-1}(P)=P$
and $\phi_3=\phi|P=I|P=L_3$. 
Now consider the box $B_\infty=\{x\in \mathbb R^2: |x_1|\leq 1, |x_2|\leq 1\}$, then
$L_1(B_\infty \cap F_1') = [(0,0),(1,1)]$ is a segment. Similarly
$L_2(B_\infty \cap F_2') = [(0,0),(0,1)]$. Moreover,
$L_3(B_\infty \cap 0^+P) = B_\infty\cap 0^+P$. The convex hull of the union of these three
polytopes is $B'={\rm co}\{(0,0),(0,1),(1,1)\}$. Hence we have
$B={\rm co}(B' \cup (-B'))= {\rm co}\{(-1,-1),(0,1),(1,1),(0,-1)\}$. The construction shows
$|\phi(q)|_\square \leq |q|_\infty + k$ for every $q$, where $|x|_\square=|Tx|_\infty$
is the polyhedral norm generated by $B$, obtained with $T=\begin{bmatrix} 2&-1\\0&1\end{bmatrix}$.

Now observe that instead of the $|\cdot|_\infty$-unit ball $B_\infty$ we can choose a larger polytope $B_\triangle$, which still satisfies
$L_i(B_\triangle \cap F_i') \subset B'$. Namely, we can take
$B_\triangle = {\rm co}(B_\infty \cup \{(2,0),(-2,0)\})$.
The polyhedral norm is $|q|_\triangle = |Sq|_\infty$ with
$S^T=\begin{bmatrix} 0&1/2&1/2\\1&1/2&-1/2\end{bmatrix}$.
By example \ref{example1}, BIBO-stability of
\begin{align*}
G: \quad \begin{split}
\dot{x} &= Ax + Bp\\
q&= Cx
\end{split}
\qquad \qquad p(t)=\phi(q(t)).
\end{align*}
can now be assessed via BIBO-stability of
\begin{align*}
\widetilde{G}: \quad \begin{split}
\dot{x} &= Ax + BT^{-1}p_1\\
q_1&= SCx\\
\end{split}
\qquad \qquad p_1(t)=T\phi(S^+ q_1(t)),
\end{align*}
where $T \circ \phi \circ S^+$ is a $|\cdot|_{\infty}$-contraction.
A sufficient condition for BIBO-stability of the loop is therefore $\|\widetilde{G}\|_{\rm pk\_gn} <1$.

We compare this to the sector characterization of the non-linearity $\phi$ from example \ref{example3}.
With $H=I_2$ we get $\phi(q)^T (\phi(q)-q)\leq 0$, 
hence
$\psi(q)=2\phi(q)-q$ is an  $L_2$-contraction.
This leads to 
\begin{align*}
\widehat{G}: \quad \begin{split}
\dot{x} &= Ax + \textstyle\frac{1}{2} BCx  + \textstyle\frac{1}{2} B\widehat{p}\\
\widehat{q}&= Cx
\end{split}
\qquad \qquad \widehat{p}(t)=\psi(\widehat{q}(t)),
\end{align*}
whence a sufficient condition for $L_2$-stability is $\|\widehat{G}\|_{L_\infty} < 1$,
and this can also be obtained from the circle criterion.
Choosing 
$$
A = \begin{bmatrix}-1&2\\0.001 & -3 \end{bmatrix},
B = \begin{bmatrix} 2&-1\\0&1 \end{bmatrix} , C=\rho \begin{bmatrix} 1&0\\1&1\end{bmatrix} 
$$
we get two curves $n_1(\rho)=\|\widetilde{G}_\rho\|_{\rm pk\_gn}$ and $n_2(\rho)=\|\widehat{G}_\rho\|_\infty$.
We have $n_1(0.499) = 0.9988$, $n_2(0.499)=1.2597$,
$n_1(0.434)=0.8687$, $n_2(0.434)=0.9995$, so in between these two values
the test $\|\widetilde{G}\|_{\rm pk\_gn}<1$ guarantees BIBO-stability, while the circle criterion fails to
prove $L_2$-stability. 
\end{example}

\begin{example}
\label{example8}
Multi-dimensional saturation is not always suitably seized by the Euclidean norm \cite{ohno}. 
For signals
$x(t)\in \mathbb R^n$ consider a convex polytope $P$ of dimension $n$ with the origin
in its interior, and let $\mu_P (x)=\inf \{\mu \geq 0: x \in \mu P\}$ be its gauge function.
Then a structured saturation operator is ${\rm sat}_P(x) = x$ for $\mu_P(x) \leq 1$,
${\rm sat}_P(x)= x/\mu_P(x)$ for $\mu_P(x) > 1$. In other words, if $x(t) \in P$, then
the signal is unaffected by saturation, but if $x(t)$ reaches the boundary of $P$, then along
every ray $\rho x$, $\rho > 0$,
the magnitude of the signal is frozen at the value  it had attained when crossing the boundary
of $P$, while the direction of the signal is unchanged. This is now a special case of
Example \ref{example1}.
\end{example}

\begin{example}
A typical case is signal clipping, where $y=\sigma(x)$ is given as
$y_i ={\rm sign}(x_i)$ if $|x_i| > 1$, $y_i=x_i$ otherwise. Here $B'=B_\infty$. Indeed, consider
for simplicity the case $n=2$. Then $\sigma$ is piecewise affine with $9$ different polytopes
$P_1,\dots,P_9$, where $P_9=B_\infty$, $P_1 = \{(x_1,x_2): -1\leq x_1 \leq 1, x_2 \geq 1\}$,
$P_2 = \{(x_1,x_2): x_1 \geq 1,x_2 \geq 1\}$, $P_3=\{(x_1,x_2): x_1 \geq 1, -1 \leq x_2 \leq 1\}$,
etc. We have $A_1(x) = (x_1,1)$, $Q_1 = [-1,1) \times \{0\}$, $C_1= \{0\} \times \mathbb R_+$,
$A_2(x) = (1,1)$, $Q_2 = \{(1,1)\}$, $C_2 = \mathbb R^+ \times \mathbb R^+$, $C_2'=\{(0,0)\}$,
etc. So only the four facets  among the 9 faces of $B_\infty$ contribute to $B'$.

This immediately applies to systems like $\dot{x} = \sigma(Ax+b)$ or $\dot{x}=\sigma(Ax +Bu)$ etc.
as for instance considered in \cite{sontag}.
\end{example}

\begin{example}
\label{example9}
The authors of \cite{convex_combin}
consider non-linear systems
$\dot{x} = \sum_{i=1}^N \mu_i(x,u) \left[A_ix + B_iu \right]$, where
$\mu_i(x,u) \geq 0$, $\sum_{i=1}^N \mu_i(x,u)=1$. This can be modeled by
an operator $p = \phi(q) := \sum_{i=1}^N \mu_i(q) q_i$ with $\mu_i(q)$ a convex combination,
so that $|\phi(q)|_1 \leq |q|_\infty$, which is a polyhedral non-linearity. 
\end{example}

\begin{example}
(Attractors, limit cycles, chaotics).
In \cite{sepulchre} the authors generate MIMO non-linearities $\Delta$
by putting 
LTI-systems $H$ in feedback with static non-linearities $\Phi$. This leads to
attractors, limit cycles, chaotic behavior, and much else. 
Some of these may be
considered special cases of  (\ref{upper_triangular}). 
The out-set is a dynamic system $z=\Delta(w)$:
\begin{align}
\label{scroll}
H: \;
\begin{split}
\dot{x} &= Ax + Bp + Bw\\
q &=Cx, \quad z=x, 
\end{split}
\qquad
p = \Phi(q)\qquad \Delta = (H,\Phi)
\end{align}
where $\Phi(0)=0$. If $|\Phi(q)|_\infty \leq r |q|_\infty + k$ asymptotically, then by Theorem \ref{small_gain},
BIBO-stability of $\Delta$ follows from $\|H\|_{\rm pk\_gn} < r^{-1}$. In particular,
if $\Phi$ has bounded range, then we can choose $r > 0$ arbitrarily,
hence $\|H\|_{\rm pk\_gn} < \infty$ gives BIBO-stability of $\Delta$.

For instance, similar to
\cite[Example 4.6]{sepulchre} we let
$\Phi(q_1,q_2,q_3)= (\phi(q_3),0,0)$ with $\phi(q) = a \tanh(k q)+\rho q$, $B=C=I_3$,  and
\begin{align*}
    A=\begin{bmatrix}
    -(\beta_1+\beta_2+\beta_3)&-(\beta_1\beta_2+\beta_1\beta_3+\beta_2\beta_3)/M&-\beta_1\beta_2\beta_3/M\\M&0&0\\0&1&0\end{bmatrix}.
\end{align*}
For $M=a=k =10$, $\rho = 0.3$, $\beta_1=2$, $\beta_2=3$, $\beta_3=5$, the non-linearity $\Delta$ has three steady
states, the unstable  $0$, and two stable attractors
$(0,0,\pm 2.963)$ (Fig. \ref{fig_sepulchre} left).
The system is globally BIBO-stable, because $\phi$ has slope $\rho$ at infinity, so
integrating $\rho q_3$ into the system $H$ gives $H_0$ in loop with a non-linearity $\Phi_0$
of bounded range. Then
$r=1/3$ gives $\|H_0\|_{\rm pk\_gn} = 2.0759 < 3=r^{-1}$.
The linearization of $\Delta$ at $0$ has system matrix $A+E_1$
with $E_1=\begin{bmatrix} 0,0, \phi'(0) ; 0,0,0;0,0,0 \end{bmatrix}$, which is unstable, but
a stable linearization can readily be obtained by
linearizing about one of the attractors, e.g. $\bar{x}=(0, 0,  2.963)$.

\begin{figure}[ht!]
\begin{center}
\includegraphics[scale=0.29]{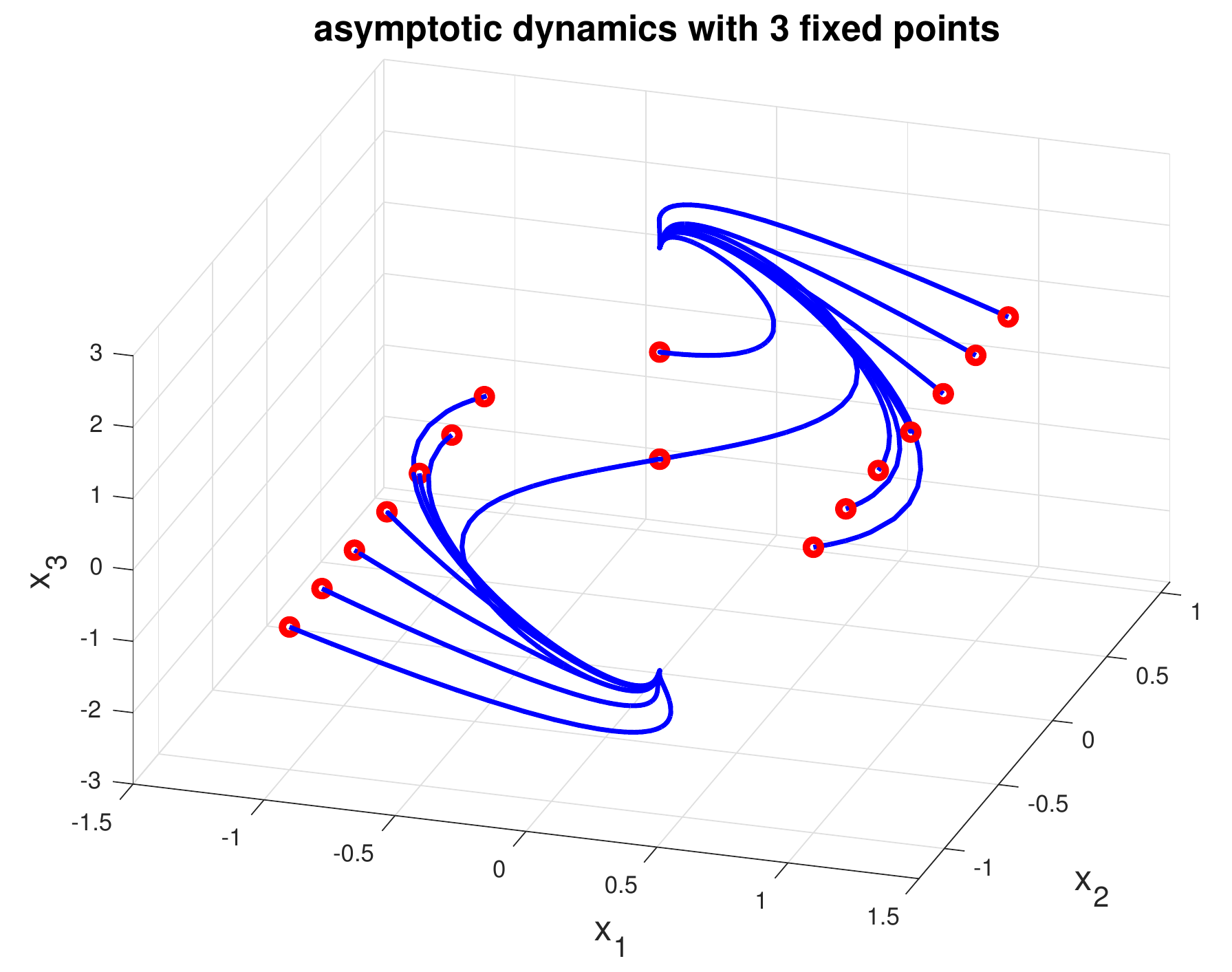}
\includegraphics[scale=0.29]{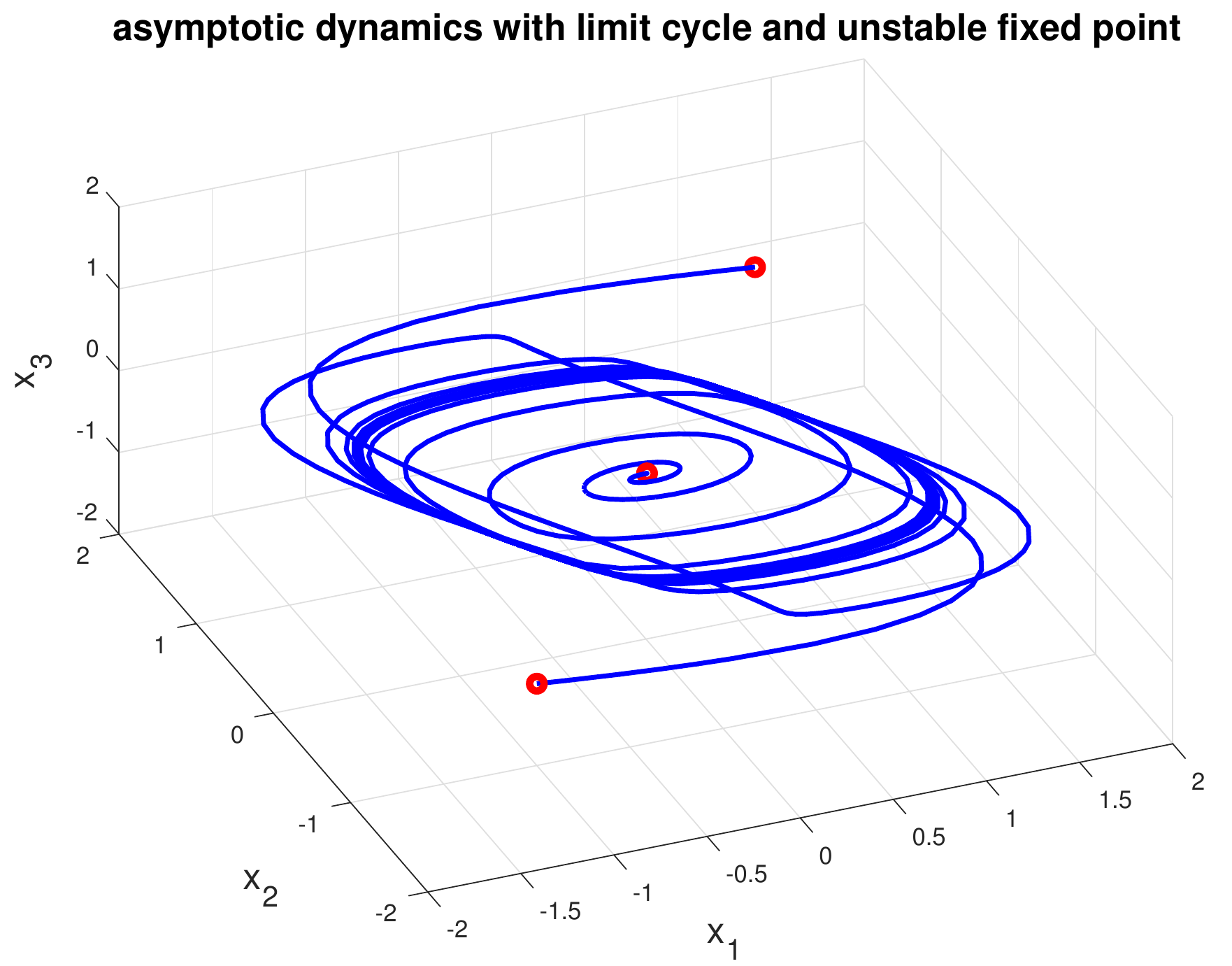}
\includegraphics[scale=0.29]{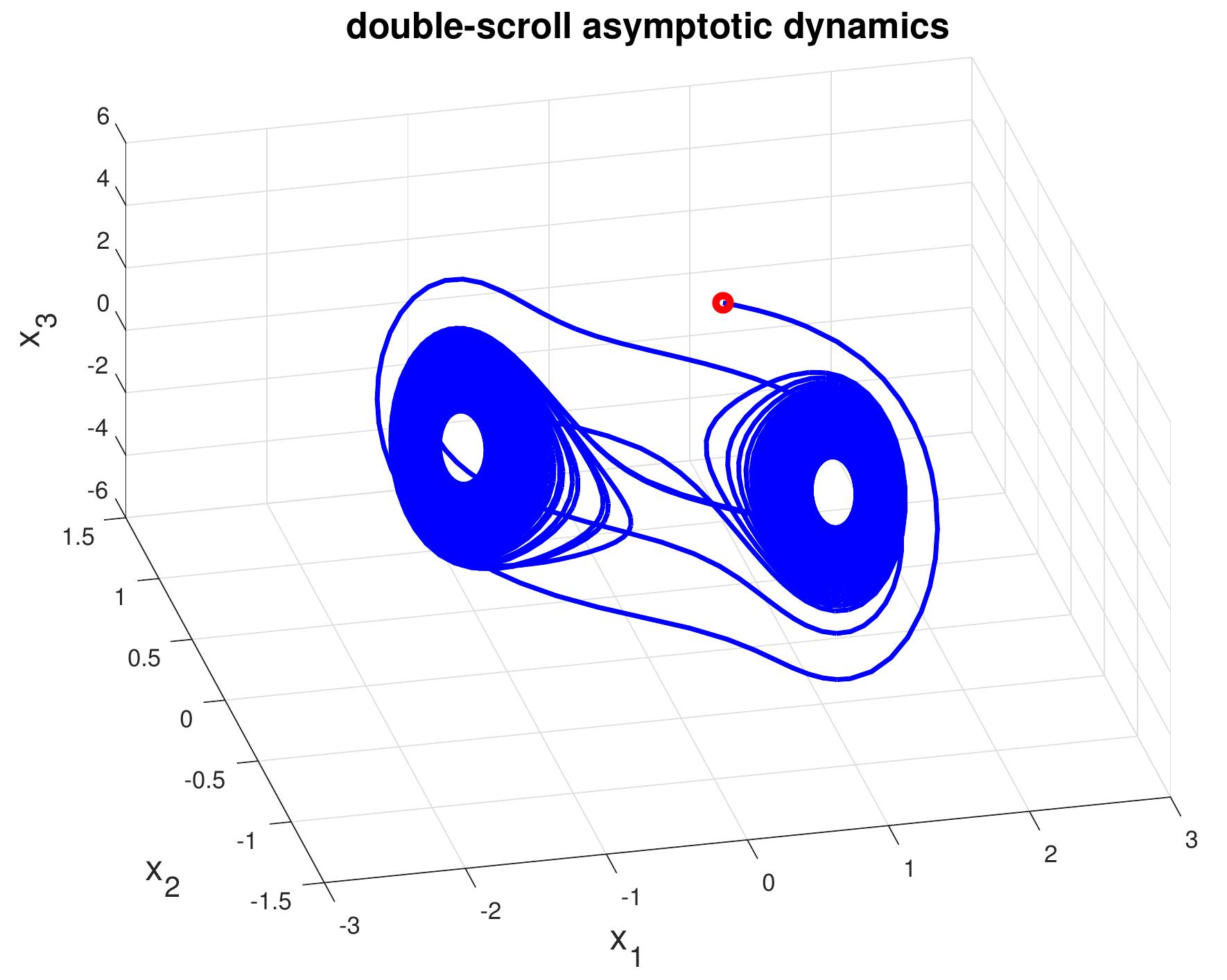}
\end{center}
\caption{Two attractors and unstable origin (left), limit cycle (middle), chaotic double scroll (right). \label{fig_sepulchre}}
\end{figure}

A different case has
$a = -10$, where
the origin is an unstable steady state and an attracting limit cycle occurs (Fig. \ref{fig_sepulchre} middle).
Global BIBO-stability of the loop follows again from $\|H\|_{\rm pk\_gn}< \infty$, but a stable
linearization now requires a modified LTI-system, where the limit cycle is subtracted from $H$
to get a stable
steady state. 
\end{example}

In those cases, where $z=\Delta(w)$ is BIBO-stable, we can consider it in loop with a
tunable LTI-system $G$ as in Fig. \ref{fig_pos_neg}. For unstable $H$ we can still
apply the results of section \ref{4.2} when $H$ is stabilized by feedback with a tunable $G$,
now considering $\Phi$ in loop with $\mathcal F_l(H,G)$:

\begin{proposition}
\label{prop2}
Let $z=\Delta(w)$ in {\rm (\ref{scroll})} be in upper feedback with an LTI-system $G$ as in {\rm (\ref{in_out})}.
Suppose $\Phi$ satisfies an asymptotic 
constraint $|\Phi(q)|_\infty \leq r|q|_\infty + k$. Then a sufficient
condition for global BIBO-stability of the loop $(G,\Delta)$
is $\|\mathcal F_l(H,G)\|_{\rm pk\_gn} < r^{-1}$.
\end{proposition}

\begin{example}(Attractors $\dots$ continued).
An interesting study in this line is Chua's circuit \cite{chua}, see \cite[5.4]{sepulchre}, where
$$A=\begin{bmatrix} -\alpha & \alpha &0\\1&-1&1\\0&-\beta&0\end{bmatrix}, B = C= I_3,
\begin{array}{l}
\Phi(q)=(\phi(q_1),0,0) \\ \vspace{.1cm} \phi(q_1)=\alpha \tanh(2q_1)+\alpha\rho q_1
\end{array}$$ 
For $\alpha = 8.3$, $\beta =16.5$, $\rho=0.25$ the double scroll attractor appears
(Fig. \ref{fig_sepulchre} right). The non-linearity
has slope $\alpha\rho$ at infinity. Therefore global BIBO-stability of $\Delta=(H,\Phi)$
follows from stability of $A+E_2$, $E_2 = \begin{bmatrix} \alpha\rho,0,0;0,0,0;0,0,0\end{bmatrix}$.

A common feature of these examples is that the sector
non-linearity invites attempting $L_2$-stability, which however fails
due to the persistence of more than one attractor in feedback. This is where global
$L_\infty$-stability is still in business.

\begin{example}
(Attractors $\dots$ continued).
 We study this phenomenon in more detail through a feedback design example.
Consider the MIMO Lur'e system
\begin{align}
\label{scroll2}
H: \;
\begin{split}
\dot{x} &= Ax + Bp + B_u u\\
q &=Cx, \quad y= C_y x, 
\end{split}
\qquad
p = \Phi(q)
\end{align}
with $\Phi:\mathbb R^3 \to \mathbb R^3$ the MIMO static non-linearity
$$
\Phi(q) := \begin{bmatrix}
\frac{q_1^2}{a_1+q_1^2} (\tanh(c_1q_1) +\rho_1 q_1))\\
\frac{q_2^2}{a_2+q_2^2} (\tanh(c_1q_2) +\rho_2 q_1)) \\
\frac{q_3^2}{a_3+q_3^2} (\tanh(c_3q_3) +\rho_3 q_3))
\end{bmatrix}\,,
$$
with $a_1 = 0.1$, $a_2 = 0.2$, $a_3 = 0.3$, $\rho_1 = 0.1$, $\rho_2=0.2$, $\rho_3 = 0.3$, $c_1=2$, $c_2 = 3$ and $c_3=4$. 
State-space data are given as
$$
A: = \begin{bmatrix}
-2 &8.8 & 0\\
1&-1&1\\
0&-15&0\\
\end{bmatrix}\,, B: = \begin{bmatrix}
5 &  0 &  0 \\
            0     &       0.1    &        0 \\
            0     &       0     &       0.3 \\
\end{bmatrix}\,, C: = I_3\,,
$$
$$
B_u := \begin{bmatrix} 1  \\
     1    \\
     1   \end{bmatrix}\, C_y := \begin{bmatrix}1   &  1   &  1  \end{bmatrix}\,.
$$
The uncontrolled linear dynamics show $2$ unstable oscillating modes $0.1422 \pm 3.0189i$. 
Simulations of the uncontrolled non-linear system are shown in Fig. \ref{2scroll} (upper line) 
with a  double-scroll regime close to the origin, and an escaping unstable spiral regime away from $0$. 
Here we regard $\Delta=(H,\Phi)$ as mapping initial conditions $x(0)$ to state $x$.

Now we investigate whether the system may be stabilized by feedback $u=Ky$  in the $L_\infty$-sense,
using the techniques in sections \ref{4.1} and  \ref{4.2}. 
As can be seen each component $\Phi_i$ belongs to the asymptotic sector  ${\bf sect}(0, \rho_i+\epsilon)$ for any $\epsilon >0$. Using section \ref{4.2}, we infer that
the closed-loop system is $L_\infty$-stable whenever $\mathcal  \|\mathcal F_l(\widetilde{H},K)\|_{\rm pk\_gn} < r^{-1}$, where  $r:=\max \left\{ \rho_1/2,\rho_2/2,\rho_3/2  \right\}$ and $\widetilde{H}$ is obtained from $H$ by centering the non-linearity. The latter amounts to shifting the $A$-matrix 
to $A +  B\Gamma C$ with $\Gamma:= {\rm diag}\left({\rho_1}/{2},{\rho_2}/{2},{\rho_3}/{2}\right)$.

Running  program (\ref{second})  over the class $\mathscr K$ of PID controllers leads to 
$$
K^*(s):= -0.796 + \frac{0.000352}{s} + \frac{0.097}{940s+1}, 
$$
with the result $\| \mathcal F_l(\widetilde{H},K^*)\|_{\rm pk\_gn} = 5.34 < r^{-1}= 6.67 $, 
affirming BIBO stability. 
Closed-loop simulations in Fig. \ref{2scroll} bottom show co-existence of $3$ stable equilibrium points at the origin (right) and away from the origin  (bottom left) with state-space coordinates $(\pm 2.98, \mp 0.0420,\mp 2.94, \mp 237.94,0)$. 

Note that in this study step $4$ of the algorithm fails, 
because $\Phi_i \in {\bf sect}(0, 1.17)$ tightly for every $i$. Therefore, to get a global
$L_2$-stability certificate it would have been required to
determine a PID controller for which the $H_{\infty}$ norm of the  
corresponding centered system was less than $1/r_0= 1.71$, and this value was
not achievable in program (\ref{first}).

\begin{figure}[ht!]
\begin{center}
\hspace{-1cm}\includegraphics[width=1.1\textwidth ]{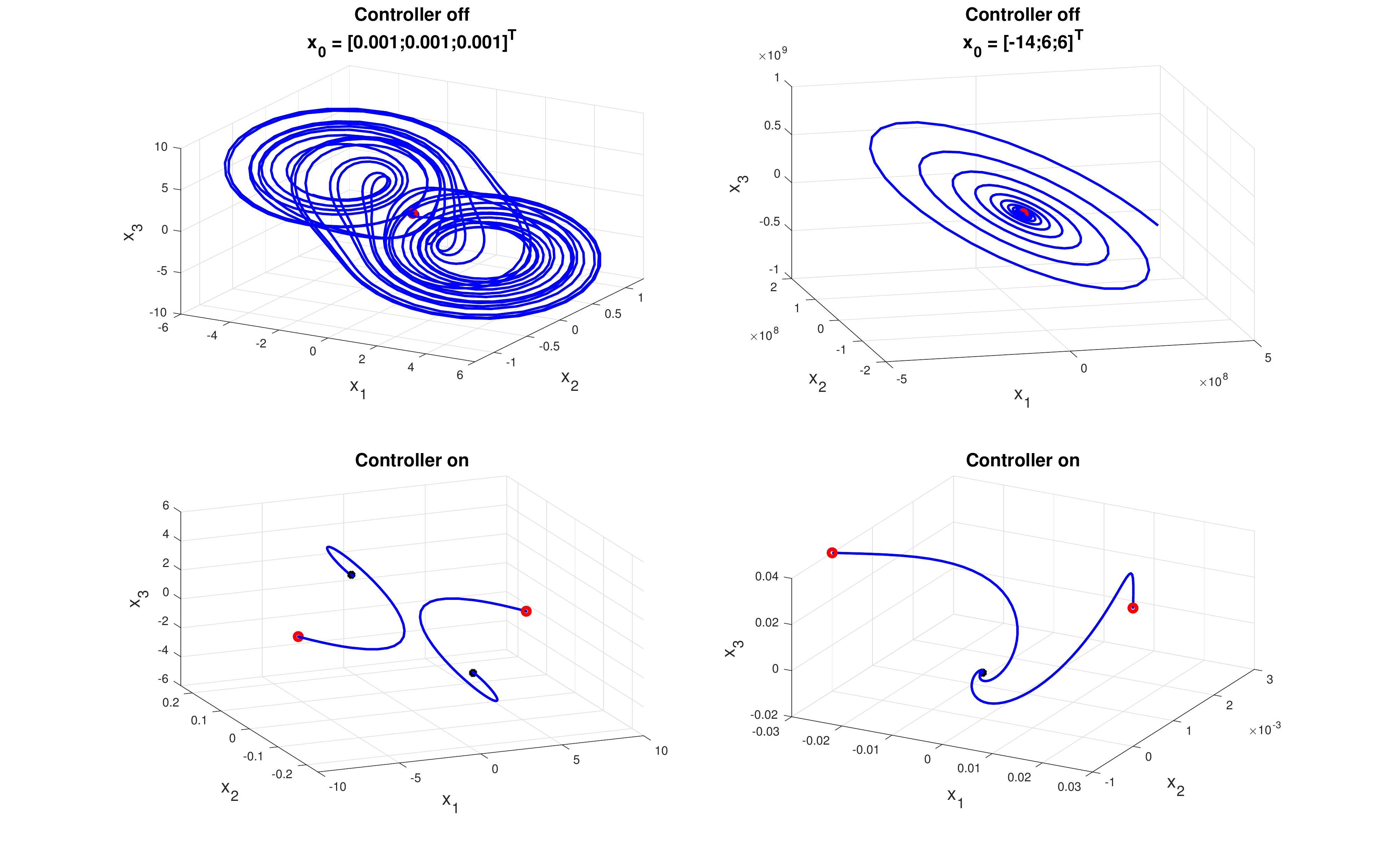}
\end{center}
\caption{Non-linearity $\Delta = (H,\Phi)$ with unstable $H(s)$ shows $2$-double scroll (top left) or diverging spiral (top right) for different initial values '$\color{red}{\bullet}$'.
Peak-gain optimization with a PID controller gives  BIBO-stability of $(K,\Delta)=((K,H),\Phi)$  with convergence to $3$ stable equilibrium points (bottom left and right). '${\color{blue}\bullet}$'.  \label{2scroll}}
\end{figure}
\end{example}

\end{example}

\section{Peak-to-peak norm}
\label{peak}

\subsection{Estimate}

It is well-known \cite{induced_norms,dahleh,desoer,boyd_barratt,boyd_doyle} that
for real-rational systems $G$ the peak-gain or peak-to-peak norm
is
\begin{equation}
    \label{pk_gn}
\|G\|_{\rm pk\_gn} = \max_{i=1,\dots,m} \sum_{j=1}^p \left( |g_{ij}^0|_1 + |d_{ij}|\right),
\end{equation}
where  $g_{ij}(t)=c_ie^{At} b_j + d_{ij}\delta(t)=g_{ij}^0(t)+d_{ij} \delta(t)$ with $g_{ij}^0\in L_1$.
A special case is the well-known expression 
$$
\normi{A}_\infty
= \max_{i=1,\dots,m} \sum_{k=1}^p | a_{ik}| = \max_{i=1,\dots,m} |{\rm row}_i(A)|_1
$$
of the maximum row-sum-norm of
$A \in \mathbb R^{m\times p}$, i.e., the induced $\ell_\infty$-$\ell_\infty$ matrix norm.

Formula
(\ref{pk_gn}) holds also for infinite dimensional BIBO-stable
systems and may be justified e.g. by the approach  \cite{unser,feichtinger}, which considers
BIBO-stable systems as all those LTI-systems
$G$, where $G(s)$ is the Laplace transform of a matrix-valued Radon measure of bounded variation.
While \cite{unser} handles the SISO case, where the 
norm is referred to as the $\mathcal M$-norm
$\|\mu\|_\mathcal M = \sup_{\phi\in \mathscr D(\mathbb R^+),\|\phi\|_\infty \leq 1} \int_{\mathbb R^+} \phi d\mu$,
one easily generalizes this to the MIMO case and obtains the formula
$$
\|G\|_{\rm pk\_gn} = \sup_{i=1,\dots,m} \sum_{k=1}^p \|\mu_{ik}\|_\mathcal M,
$$
which contains (\ref{pk_gn}) as a special case. In particular, it was possible to
apply it in the slipstick study, because the impulse response was an element of $L_1$.

Estimate (\ref{chain}) is mentioned in \cite{dahleh} with non-specified constants,
and the SISO case $p=m=1$ is proved in \cite[Thm.]{boyd_doyle} for  discrete SISO systems,
and in \cite[pp. 11-12]{mil} for continuous SISO systems, where in the latter reference the idea of proof is attributed to I. Gohberg. The left hand estimate in (\ref{chain}) holds also for infinite dimensional systems, e.g. those 
where $G(s)$ is the Laplace transform of a Radon measure of bounded variation,
while the right hand estimate is true for finite-dimensional $G$.
For strictly proper
systems, the estimate $\|G\|_{\rm pk\_gn} \leq 2np^{1/2}\|G\|_\infty$ is valid.
Details on computing these estimates will be published elsewhere.

\subsection{Implementation}

\label{new_section}
Stand alone
computation of the peak-gain  norm with high precision has been addressed in the literature
\cite{rutland,linnemann,boyd_doyle,induced_norms,dahleh,dahleh_pears,diaz}. For optimization,
due to non-smoothness of both norms in (\ref{second}), we need to supply subgradients
of closed-loop integral functionals 
$\phi_{ij}:K \to \int_0^\infty |c_i(K)e^{A(K)t}b_j(K)| dt$, those for the $H_\infty$-norm being well-known
\cite{AN06a}. Putting
$f(K,t) = c(K)e^{A(K)t}b(K)$ and $F(K)=\int_0^\infty |f(K,t)| dt$ for the generic terms, we need 
partial derivatives $\partial f(K,t)/\partial K_{\mu\nu}$, where $K_{\mu\nu}$ are the 
controller gains, which depend in turn on the tunable parameters ${\bf x}$
over which we ultimately optimize. Since dependence $K_{\mu \nu}({\bf x})$ on ${\bf x}$
and $f(K,t)$ on $K$ is differentiable, non-smoothness occurs only when the absolute value
$|f(\cdot,t)|$ is formed, and ultimately via the finite maximum over rows in (\ref{pk_gn}).
Subgradients $g(K,t) \in \partial |\cdot| \circ f(K,t)$
are obtained as $g_{\mu \nu}(K,t) = \partial f(K,t)/\partial K_{\mu\nu} {\rm sign} f(K,t)$ for 
$f(K,t) \not=0$, while those $(K,t)$ where $f(K,t)=0$ give
the full set of subgradients $g_{\mu\nu}(K,t) \in \partial f(K,t)/\partial K_{\mu\nu}\cdot [-1,1]$.
Partial derivatives $\partial f/\partial K_{\mu\nu}$
are obtained via algorithmic differentiation \cite{alberto2}. Finally,
subgradients $G \in \partial F(K)$ of integral functionals
are by regularity simply integrals of pointwise subgradients
$G_{\mu\nu}(K) = \int_0^\infty  \partial f(K,t)/\partial K_{\mu\nu} {\rm sign} f(K,t) dt$ \cite{giner,noll1,noll2}.

For
mixed programs like (\ref{second}) it is possible to use a progress function approach
as in \cite{simoes1,simoes2,simoes3,dao1,dao2,dao3,aude,gabarrou}. 
Here we rather follow the line  \cite{an_disk:05,laleh1}
suited for norm functionals, where an iteratively re-weighted
maximum of several norms is minimized; cf. \cite{encyc} for an overview, where in particular the
concept of hard and soft constraints is addressed. This approach is also used in 
the {\tt systune} function \cite{gahinet2012frequency,apkarian2014multi,CST2020b} based on \cite{an_disk:05}, and has been used in our experiments.
For convergence issues of bundle and bundle trust-region techniques
we refer to \cite{noll3,ANR:2015}.

According to the line in \cite{an_disk:05}
program (\ref{second}) is addressed by minimizing a maximum
$$
\min_{\x \in \mathbb R^n} \max\left\{ \alpha \|T_{wz}(G,K(\x))\|_\infty,\beta \|T_{pq}(G_\psi,K(\x))\|_{\rm pk\_gn}\right\},
$$
where the weights $\alpha,\beta$ are updated iteratively until the constraint of
(\ref{second}) is satisfied, from where on the objective is reduced. Here
$K(\x)$ expresses dependence of $K$ on the tunable parameters $\x$. The first 
term  splits into a semi-infinite maximum
$$
\|T_{wz}(G,K(\x))\|_\infty=\max_{\omega \in [0,\infty]} \overline{\sigma}(T_{wz}(j\omega,G,K(\x))),
$$
whereas the second term, after time-domain discretization,  becomes a finite maximum
$$
\|T_{pq}(G_\psi,K(\x))\|_{\rm pk\_gn}=
\max_{i=1,\dots,m}\sum_{j=1}^p \left(\sum_{t\in T}\left| c_i(\x)e^{A(\x)t}b_j(\x)\right| + |d_{ij}(\x)|\right),
$$
with $c_i(\x)e^{A(\x)t}b_j(\x) + d_{ij}(\x)\delta(t)$ the closed loop impulse response of the entry $(i,j)$ of the channel
$T_{pq}(G_\psi,K(\x))$. For discretization we have used the method of \cite{rutland}, which is 
readily extended to the 
MIMO case.

It is helpful to update the weights $\alpha,\beta$ in such a way that at the current iterate 
$\x$ the two branches are at least nearly active. Selecting a set of active and near active frequencies for the first 
objective is explained in \cite[sect. 4.4]{an_disk:05}, and we proceed analogously
for the second branch.  This strategy to include near-active branches into local models
has turned out highly effective, as it
avoids stalling at non-optimal points.

\section{Conclusion}
We have presented a method for stabilization and performance optimization
of non-linear controlled systems, where the non-linearity satisfies a sector constraint
asymptotically. This leads to global closed-loop BIBO-stability in tandem with local exponential stability
in  situations where global closed-loop $L_2$-stability fails, either due to
exceedingly large sectors, or more principally, due to persistence of 
several attracting regimes in closed loop. The new approach requires solving
a mixed $L_1/H_\infty$-synthesis program, and uses properties of the $L_1$- or
peak-gain system norm.


\end{document}